\documentclass{article}
\usepackage[lang=american]{ems-journal}
\usepackage{enumitem}
\usepackage{mathtools}
\usepackage{mathrsfs}
\DeclareFontFamily{U}{rsfs}{\skewchar\font127 }
\DeclareFontShape{U}{rsfs}{m}{n}{%
  <-6.5> rsfs5 <6.5-8> rsfs7 <8-> rsfs10}{}
\usepackage{microtype}
\usepackage{xurl}
\numberwithin{equation}{section}
\theoremstyle{plain}
\newtheorem{theorem}{Theorem}[section]
\newtheorem{proposition}[theorem]{Proposition}
\newtheorem{lemma}[theorem]{Lemma}
\newtheorem{corollary}[theorem]{Corollary}
\theoremstyle{definition}

\theoremstyle{definition}

\newtheorem{rem}[theorem]{Remark}
\newcommand{\bE}{\mathbb E}
\newcommand{\bC}{\mathbb C}
\newcommand{\bD}{\mathbb D}
\newcommand{\cH}{H}
\newcommand{\1}{\mathbf 1}
\newcommand{\pd}{\partial}
\newcommand{\eps}{\varepsilon}
\newcommand{\om}{\omega}
\newcommand{\R}{\mathbb R}
\newcommand{\C}{\mathbb C}
\newcommand{\T}{\mathbb T}

\newcommand{\KD}{K_D}
\newcommand{\cD}{a}
\newcommand{\Gcat}{\beta(2)}
\newcommand{\Bell}{\mathscr B}

\newcommand{\Hilb}{H}

\allowdisplaybreaks[1]

\hypersetup{pdftitle={On the sharp weak-type inequality for Riesz transforms},pdfauthor={Rodrigo Banuelos, M\"ahdi Hormozi, Alexander Volberg}}
\ifdefined\XeTeXversion
  \catcode"2019=\active
  \protected\def’{'}
  \catcode"2013=\active
  \protected\def–{--}
\fi
\begin{document}

\title{On the sharp weak-type inequality for Riesz transforms}
\titlemark{On the sharp weak-type inequality for Riesz transforms}

\emsauthor{1}{
  \givenname{Rodrigo}
  \surname{Ba\~nuelos}}{R.~Ba\~nuelos}
\Emsaffil{1}{
  \department{Department of Mathematics}
  \organisation{Purdue University}
  \zip{47907}
  \city{West Lafayette, IN}
  \country{USA}
  \affemail{banuelos@purdue.edu}}

\emsauthor{2}{
  \givenname{M\"ahdi}
  \surname{Hormozi}}{M.~Hormozi}
\Emsaffil{2}{
  \organisation{Beijing Institute of Mathematical Sciences and Applications}
  \zip{101408}
  \city{Beijing}
  \country{China}
  \affemail{hormozi@bimsa.cn}}

\emsauthor{3}{
  \givenname{Alexander}
  \surname{Volberg}}{A.~Volberg}
\Emsaffil{3}{
  \department{Department of Mathematics}
  \organisation{Michigan State University}
  \address{619 Red Cedar Road}
  \zip{48824}
  \city{East Lansing, MI}
  \country{USA}
  \affemail{volberg@math.msu.edu}}

\begin{abstract}
We prove that the weak-type $(1,1)$ norm of each classical Riesz
transform $R_j$, $j=1,\dots,n$, on real-valued inputs on $\R^n$ and
$\T^n$ is $\KD$. The constant $\KD=\pi^2/(8\Gcat)$, where $\Gcat$
is Catalan's constant, was obtained by Davis \cite{Davis1974} for the
Hilbert transform on $\R$ and the conjugate-function operator on the
circle $\T$. The Burkholder--Bellman function argument also applies
to coordinate Riesz transforms of product generators. If $-L$
generates a symmetric sub-Markov semigroup on a nonzero standard
$\sigma$-finite measure space $Y$, then
$-\partial_x(-\partial_x^2+L)^{-1/2}$ has the same exact norm on both
$\mathbb R\times Y$ and $\mathbb T\times Y$. The derivative acts in
the ordinary line or circle factor, while $L$ governs the independent
transverse variables. An explicit integral inequality strengthens
the weak-type estimate and bounds the truncated second moment of the
transform. The general theorem applies to Gaussian, Dunkl, Laguerre, Jacobi,
Schr\"odinger, and Laplace--Beltrami transverse generators,
finite reversible chains, and subordinated semigroups.
\end{abstract}

\classification[42B15, 47D07]{42B20}


\makeatletter
\let\ps@titlepage\ps@empty
\makeatother
\maketitle

\section{Introduction}\label{sec:introduction}\label{hv-product-section}

The weak-type (1,1) sharp constant problem for the Riesz transforms has its origins in
the classical inequalities of M.~Riesz and Kolmogorov for conjugate
functions. For $j=1,\dots,n$, we write
$R_j=-\partial_j(-\Delta)^{-1/2}$, with Fourier multiplier
$-i\xi_j/|\xi|$ on $\R^n$ and $-im_j/|m|$ on $\T^n$, zero at the
origin. Our Euclidean Fourier convention is
$\widehat f(\xi)=\int e^{-ix\cdot\xi}f(x)\,dx$.
On $\R^n$, these are the Calder\'on--Zygmund singular integrals
\begin{equation*}
 R_jf(x)=c_n\,\operatorname{p.v.}\int_{\R^n}
 \frac{y_j}{|y|^{n+1}}f(x-y)\,dy,
 \qquad c_n=\pi^{-(n+1)/2}\Gamma\bigl((n+1)/2\bigr).
\end{equation*}
For $n=1$, these are the Hilbert transform $H$ on $\R$ and
the conjugate-function operator on $\T$.

For a scalar or Hilbert-space-valued function $v$,
\[
 \|v\|_{1,\infty}=\sup_{\lambda>0}\lambda\,\mu\{|v|>\lambda\},
\]
where $|v|$ is the absolute value or Hilbert norm.
The sharp endpoint estimates proved below concern real-valued inputs.

Davis \cite{Davis1974} determined the optimal constant in
Kolmogorov's inequality on both $\R$ and $\T$:
\[
 \KD=\frac{\pi^2}{8\Gcat},
 \qquad \Gcat=\sum_{k\ge0}\frac{(-1)^k}{(2k+1)^2}.
\]

\begin{theorem}\label{thm:classical}
Let $n\ge1$, let $1\le j\le n$, and let $\Omega=\R^n$ or $\T^n$.
For every real-valued $f\in L^1(\Omega)$,
\begin{equation}\label{eq:classical-weak}
 \|R_jf\|_{1,\infty}\le\KD\|f\|_1,
 \qquad \KD=\frac{\pi^2}{8\Gcat}.
\end{equation}
The constant $\KD$ is best possible for every $n$ and $j$.
\end{theorem}

The proof gives the following more general result for coordinate
Riesz transforms of product generators.

Let $(Y,\nu)$ be a nonzero standard $\sigma$-finite measure space,
and let $L\ge0$ be self-adjoint on $L^2(Y,\nu)$. We assume that
\[
 P_t=e^{-tL}
\]
is a symmetric sub-Markov semigroup: it preserves positivity and
satisfies $P_t1\le1$ in its consistent $L^\infty$ realization.
The conservative case has $P_t1=1$; the inequality also permits killing.

On either $X=\T\times Y$ or $X=\R\times Y$, with product measure
$\mu$, put
\begin{equation*}
 A=(-\partial_x^2+L)^{1/2},\qquad R=R_L=-\partial_xA^{-1}.
\end{equation*}
The circle $\T=\R/(2\pi\mathbb Z)$ has normalized Haar measure
$dm(x)=dx/(2\pi)$, and the line has Lebesgue measure.
We call $x$ the horizontal variable and the variables in $Y$
transverse. No differentiation in $Y$ is required.

The definition of $R$ is spectral. On the circle,
\begin{equation*}
 \widehat{Rf}(m)=-im(m^2+L)^{-1/2}\widehat f(m)
 \quad(m\ne0),\qquad \widehat{Rf}(0)=0.
\end{equation*}
On the line, replace $m$ by $\xi$, using
$\widehat f(\xi)=\int_\R e^{-ix\xi}f(x)\,dx$.
In particular, $R^*=-R$ and $\|R\|_{2\to2}\le1$.
The product semigroup is
\[
 e^{-tA^2}=e^{t\partial_x^2}\otimes P_t.
\]

Throughout the paper, $a=\KD^{-1}=8\Gcat/\pi^2$.

\begin{theorem}\label{thm:horizontal}
Under the preceding assumptions, $R$ has a unique extension from
real $L^1(X)\cap L^2(X)$ to real $L^1(X)$ that is continuous from
$L^1$ convergence to convergence in measure. For every real
$f\in L^1(X)$ and every $\lambda>0$,
\begin{equation}\label{eq:horizontal-weak}
 \lambda\,\mu\{|Rf|>\lambda\}\le\KD\|f\|_1.
\end{equation}
The weak-type constant is sharp:
\begin{equation*}
 \|R\|_{L^1(X;\R)\to L^{1,\infty}(X)}
 =\KD=\frac1{\cD}.
\end{equation*}
\end{theorem}

Taking $L=-\Delta_y$ gives Theorem~\ref{thm:classical}.
For $n=1$, take a one-point transverse space and $L=0$.
Theorem~\ref{thm:horizontal} does not require the semigroup
to preserve constants.
Examples include
horizontal products with Gaussian, Dunkl, Laguerre, Jacobi,
Schr\"odinger and manifold semigroups, reversible chains, and
Bernstein subordinations. In each case the numerator is the
ordinary derivative in the added line or circle variable.

An integral refinement of the weak-type estimate is stated in
Theorem~\ref{thm:horizontal-integral}, after the Davis function
has been constructed.

\subsection{Historical background}\label{sec:organization}

M.~Riesz \cite{Riesz} proved the $L^p$ boundedness of the
conjugate-function operator for $1<p<\infty$; the corresponding
theorem holds for the Hilbert transform on $\R$.
At $p=1$, Kolmogorov \cite{Kol} established the weak-type
inequality for conjugate functions. Its real-line form is
\begin{equation*}
 \|Hf\|_{1,\infty}\le C_1\|f\|_1.
\end{equation*}
The corresponding strong and weak estimates for
Riesz transforms in higher dimensions, with constants depending on $n$,
follow from the work of Calder\'on and Zygmund \cite{CZ} on a more
general class of singular integral operators.

These results led to the question of the optimal constants.
Gohberg and Krupnik \cite{GohbergKrupnik1968} determined the
$L^p$ norm of the Hilbert transform for $p=2^k$, $k\ge1$, using
Cotlar's identity. Pichorides \cite{Pichorides1972} proved that
\begin{equation*}
 \|H\|_{L^p(\R;\R)\to L^p(\R)}=
 \begin{cases}
  \tan\left(\dfrac{\pi}{2p}\right),&1<p\le2,\\[3pt]
  \cot\left(\dfrac{\pi}{2p}\right),&2\le p<\infty.
 \end{cases}
\end{equation*}
The same constants hold for the conjugate-function operator on
$\T$. Iwaniec and Martin \cite{IwMa} obtained the same sharp
$L^p$ norm for each Riesz transform on $\R^n$ by the Calder\'on--Zygmund method of 
rotations which reduces the problem to the case of the Hilbert transform.  An alternative probabilistic proof was given by Ba\~nuelos and Wang
\cite{BanuelosWang1995}
based on their inequalities for orthogonal martingales and the Gundy--Varopoulos \cite{GundyVaropoulos1979} representation of the Riesz transforms.

Dimension-free $L^p$ estimates for Riesz transforms have also been
studied in Gaussian spaces, Lie groups, and Riemannian manifolds. 
Important approaches combine the probabilistic representation of
Gundy and Varopoulos \cite{GundyVaropoulos1979} with 
martingale inequalities \cite{BanuelosWang1995,Burkholder1984,BanuelosOsekowski2015}, or use
Bellman functions and bilinear embedding estimates, following the work of
Nazarov, Treil, and Volberg \cite{NazarovTreilVolberg1999} and
Nazarov and Volberg \cite{VolbergNazarov2004}. Under suitable geometric
hypotheses, these methods give bounds depending only on $p$, often
with sharp or nearly sharp dependence on $p$. Further examples include the multiplier theorems of Ba\~nuelos, Baudoin, Chen, and Sire
\cite{BanuelosBaudoinChenSire2021}, and the Bellman-function estimates of
Carbonaro and Dragi\v{c}evi\'c \cite{CarbonaroDragicevic2013}.
For recent dimensionless weighted estimates for the Bakry--Riesz
vector obtained by continuous sparse domination, see Domelevo,
Petermichl, and \v{S}kreb \cite{DomelevoPetermichlSkreb2025}.

There are other proofs of Davis's theorem, including the complex-analytic
proof due to Baernstein \cite{Baernstein1978}. Janakiraman
\cite{JanakiramanWeak2004} subsequently determined the sharp
weak-type $(p,p)$ constants, $1\le p\le2$, for orthogonal harmonic
functions and martingales, extending Davis's result and yielding the
corresponding sharp inequalities for the Hilbert transform.

In contrast to the sharp $L^p$ theory, the corresponding weak-type
$(1,1)$ bound with constant $\KD$ for individual Riesz transforms in
dimensions $n\ge2$ had remained open. The natural conjecture was that
the weak-type $(1,1)$ norm of each Riesz transform on real-valued
inputs is also $\KD$, in every dimension. This question was asked explicitly in
\cite[Problem~7]{Banuelos2010}. The method of rotations and the
martingale arguments used in the sharp $L^p$ theory do not directly
yield the sharp weak-type $(1,1)$ estimate.

The more basic question of whether the Riesz transforms satisfy
weak-type $(1,1)$ inequalities with constants independent of dimension
was raised by Stein in his 1986 ICM address
\cite[p.~203, Problem~(b)]{SteinICM1986}. The identification of the
optimal constant is a further question.
Janakiraman obtained a weak-type bound growing logarithmically
with dimension \cite{Janakiraman2004}. Spector and Stockdale
reduced the dimension-free question to finite positive sums of
point masses \cite{SpectorStockdale2021}. More recently, Ouyang,
Spector, and Stockdale obtained the dimension-free upper bound $2$
for the full Euclidean vector Riesz transform \cite{OSS2026}.
We prove that the exact constant for each coordinate transform
on real-valued inputs, on both $\R^n$ and $\T^n$, is $\KD$.
The sharp vector constant remains a separate question which is not addressed  in this paper. 

Os\k{e}kowski
\cite[Theorem~1.1]{OsekowskiLog2012} obtained sharp $L\log L$
inequalities for Riesz transforms, with optimal constants independent
of dimension. These control the integral of $|R_jf|$ over sets of finite
measure under logarithmic integrability assumptions on $f$.
His proof uses inequalities for orthogonal, differentially subordinate
martingales, the Gundy--Varopoulos representation, and duality.
The convex integral estimates involved are preserved under conditional
expectation. This mechanism does not directly transfer sharp weak-type
$(1,1)$ estimates, since conditional expectation is not contractive in
$L^{1,\infty}$. Thus these logarithmic inequalities address a different
endpoint question from the one considered here.

Our proof is based on the Burkholder--Bellman function technique.
Burkholder's boundary-value method \cite{Burkholder1984} led to
sharp martingale inequalities, and Bellman functions are now
widely used in sharp inequalities in harmonic analysis. For applications
of this method, we refer the reader to
\cite{NazarovTreilVolberg1999,Banuelos2010,Osekowski2012,OsekowskiSurvey2013,VasyuninVolberg2020}
and references therein.
Here we use Davis's strip function directly at the operator level.
The search for the sharp constant also gives an integral inequality
that controls the truncated second moment of the transform.

Section~\ref{sec:davis} recalls Davis's strip function and gives the
Bellman proof of the one-dimensional inequality, then constructs
$\Phi$ and $\Bell$. The stronger estimate is
Theorem~\ref{thm:horizontal-integral}; inequality
\eqref{eq:main-truncated} strengthens \eqref{eq:classical-weak} and
\eqref{eq:horizontal-weak} by controlling the truncated second moment.
Section~\ref{sec:finite} proves Theorem~\ref{fp-finite} for finite
transverse spaces. The strip estimate in Subsection~\ref{fp-finite-section}
is extended by regularization, maximization, and confinement in
Subsections~\ref{fp-remove-section} and~\ref{fp-contact-section}.
Subsection~\ref{tr-section} passes to general transverse spaces by approximation and
extends the estimates to $L^1$: Proposition~\ref{tr-upper} proves the
upper estimates in Theorems~\ref{thm:horizontal} and
\ref{thm:horizontal-integral}, and Proposition~\ref{tr-sharpness}
gives the matching lower bound by horizontal oscillation at high frequencies. This also
proves Theorem~\ref{thm:classical}. Section~\ref{ex-section} gives
examples of transverse generators.

\subsection*{Acknowledgments} R. Bañuelos thanks his colleague Burgess Davis for many conversations over the past forty years about his theorem and its possible extension to $\R^n$.

\subsection*{Contributions and use of AI}\label{app:contributions} The authors conceived the mathematical approach of this paper, including the use of the Burkholder–Bellman function method, directed its development and many computations were done without AI assitance. OpenAI’s ChatGPT and Codex were used extensively for mathematical discussion, exploration and refinement of proofs, calculations, literature searches, and manuscript preparation. The authors critically assessed and revised the AI-generated contributions and take full responsibility for the mathematical content and presentation.

\section{Preliminaries}\label{sec:davis}

\subsection{Davis's strip function}
\label{Davis}\label{fp-strip-section}\label{hv-davis-section}\label{Bell1}\label{dim1}

Davis's probabilistic proof of the weak $(1,1)$ inequality \cite{Davis1974}
uses L\'evy's conformal invariance property for Brownian motion.
It rests on two inequalities:
\begin{equation}
\label{Dh}
 \bE_{r+is}|X_\tau | \ge \bE_{is}|X_\tau|;
\end{equation}
 \begin{equation}
 \label{Dv}
\bE_{is} |X_\tau| \le \bE_0|X_\tau|\,.
 \end{equation}
Here $r+ is$ lies in the horizontal strip $\Pi=\{z:|\Im z|<1\}$, $Z_t= X_t+ iY_t$ is a standard Brownian motion, $\tau$ is its first exit time from the strip, and 
the subscript indicates the point where we start the Brownian motion.

He also computed $\bE_0|X_\tau| =\frac{8\Gcat}{\pi^2}$, where 
\[
\bE_0 |X_\tau|=\frac12\int_\R |t|\operatorname{sech}\left(\frac{\pi t}{2}\right)dt
=\frac{4}{\pi^2}\int_0^\infty s\operatorname{sech}s\,ds.
\]
Using $\operatorname{sech}s=2\sum_{n\ge0}(-1)^n e^{-(2n+1)s}$,
\begin{equation}\label{eq:a}
a:=\bE_0 |X_\tau|=\frac8{\pi^2}\sum_{n=0}^\infty\frac{(-1)^n}{(2n+1)^2}
=\frac{8\Gcat}{\pi^2}.
\end{equation}
Davis's argument shows that the weak norm of the Hilbert transform $\cH: L^1(\T)\to L^{1, \infty}(\T)$ is exactly $\frac{\pi^2}{8\Gcat}$.

\medskip

We prove that the same constant is the weak-type norm of the horizontal Riesz transform $R$, independently of dimension. We begin with a Bellman function proof of Davis's result.

Consider the following obstacle problem on the plane $z=x+iy$: given an obstacle 
\[
O(x, y) =\1_{|y|\ge 1} - C|x|
\]
find the smallest $C$ for which a superharmonic majorant with $\mathcal G(x,0)\le0$ exists. The value is $1/a$, and we construct the majorant below.

The solution of the problem can be represented as
\[
\mathcal G(z) =\sup_\tau \bE_z[O(Z_\tau)],\quad z=x+iy \in \Pi,
\]
where $\tau$ runs over all bounded stopping times of planar Brownian motion.
The stopping-time representation guides the construction, as in
\cite{NTV} and \cite{VasyuninVolberg2020}. It follows that
\[
\Delta \mathcal G\le 0  \quad \text{superharmonicity};
\]
and that on an open set $\Omega$ where $\mathcal G>O$, the function $\mathcal G$ is harmonic.
This suggests the candidate 
\begin{equation*}
\mathcal G(x, y) =\begin{cases} 1-C|x|, \,\, |y|\ge 1;
\\
\text{harmonic extension inside the strip }\Pi\text{ of }1-C|x|,\quad |y|<1;
\end{cases}
\end{equation*}
We need this function to be
superharmonic on the plane and satisfy
 \begin{equation}
 \label{inside}
\mathcal G(x, y) \ge -C|x| \quad |y|<1\,.
\end{equation}

Let $u(x, y)$ be the harmonic extension of the boundary function $b(x, y)=|x|$ for $x+iy \in \pd \Pi$. 
At the center of the strip,
\begin{equation*}
u(0,0) = \bE_0 |X_\tau|= a= \frac{8\Gcat}{\pi^2}\,.
\end{equation*}
Put
\[
\mathcal G(x, y) =  1-\frac1a u(x,y)
\]
and check that \eqref{inside} is satisfied with $C=1/a$, that is $a-u(x, y)\ge -|x|$. 
The Poisson kernel $P_y$ is a symmetric probability density. Consequently
\[
u(x,y)-|x|=\int_\R(|x+t|-|x|)P_y(t)\,dt.
\]
For $x>0$, its derivative is
\[
\frac{\partial}{\partial x}(u(x,y)-x)
=\int_\R\operatorname{sgn}(x+t)P_y(t)\,dt-1\le0.
\]
By evenness, the maximum over $x$ occurs at $x=0$, so
\[
u(x,y)-|x|\le u(0,y) \le u(0,0)=a. 
\]
The last inequality is \eqref{Dv}.
\medskip

We check that $\mathcal G$ is superharmonic in $\bC$. It is harmonic inside the strip $\Pi$ and, for $|y|>1$, it is superharmonic because it equals $1-\frac1a |x|$.

It remains to check the sign of the normal derivative on $|y|=1$.

Let $\Pi=\{(x,y)\in\mathbb R^2:|y|<1\}$ and let $u$ be the harmonic extension of $|x|$ from the two horizontal sides of $\Pi$:
\begin{equation*}
\Delta u=0\quad\text{in }\Pi,\qquad u(x,\pm1)=|x|.
\end{equation*}
Extend $u$ by $|x|$ outside $\Pi$. We check the normal derivative of $u$ and the superharmonicity of the globally defined pasted function
\begin{equation}
\label{eq:paste}
\mathcal G_C(x,y)=\begin{cases}1-Cu(x,y),&|y|<1,\\1-C|x|,&|y|\ge1,\end{cases}\qquad C>0.
\end{equation}
The two definitions agree continuously on $y=\pm1$.

The function $|x|$ is subharmonic on the plane, with distributional Laplacian $2\delta_{\{x=0\}}$. Its harmonic extension satisfies $u\ge |x|$ in $\Pi$, and in fact $u>|x|$ at every interior point. 
For each fixed $x_0\ne0$, the difference $v=u-|x|$ is harmonic in a half-disc centered at $(x_0,1)$ lying in the strip, positive just inside $\Pi$, and vanishes on the boundary. The Hopf boundary-point lemma therefore gives
\begin{equation*}
u_y(x_0,1^-)<0= u_y(x_0,1^+)\qquad(x_0\ne0),
\end{equation*}
which means that $\Delta u\ge 0$ on $y=1$ in the distributional sense away from $x=0$.
Similarly $u_y(x_0,-1^+)>0$ for $x_0\ne0$, and the jump at $y=-1$ is again nonnegative.
At $x_0=0$, first convolve the boundary data $|x|$ with a smooth nonnegative probability kernel in $x$. The same maximum-principle argument gives nonnegative interface jumps for these smooth convex data. Letting the kernel tend to a point mass preserves the distributional inequality.
Hence the subharmonicity of the pasted $u$ on $\bC$ is proved, and $\mathcal G=1-u/a$ is superharmonic.

\medskip

Our function $\mathcal G=1-\frac{u(x, y)}{u(0,0)}$ is superharmonic on $\bC$ and majorizes the obstacle $O$. The Bellman function proof that the weak $(1,1)$ norm of the Hilbert transform on $\T$ is bounded by $1/a= 1/u(0,0)=\pi^2/(8\Gcat)$ also uses the following property of $\mathcal G$, which follows from \eqref{Dh}:
\begin{equation}
\label{H00}
\mathcal G(x, 0) =1-\frac{u(x,0)}{u(0,0)} \le 1-\frac{u(0,0)}{u(0,0)} =\mathcal G(0,0) =0\,.
\end{equation}

We consider $g+i\cH g=: g+ih$ on the circle for a smooth real $g$, extend $g+ih$ analytically to the disk with $h(0)=0$, and use normalized measure to write
\[
\int_\T O(g+ih) \le \int_\T \mathcal G(g+ih)\le \mathcal G(g(0), h(0)) = \mathcal G(g(0), 0) \le \mathcal G(0, 0) =0\,.
\]
The left hand side is $|\{ e^{ix}: |h|\ge 1\}| - \frac1a \int_\T |g|$. Thus the norm is bounded by $\frac{\pi^2}{8\Gcat}$.
\begin{equation}
\label{cH}
|\{ e^{ix}\in \T: |\cH g|\ge \lambda\}|\le \frac{\pi^2}{8\Gcat}\frac{\|g\|_1}{\lambda},\qquad\lambda>0\,.
\end{equation}

\subsection{The explicit Bellman function}
\label{sharper}\label{hv-stronger-section}
The constant above is sharp, but \eqref{cH} admits a refinement.
The refinement and its extension to several variables use an explicit formula for $u$.

Let $\Pi=\{z=x+iy:|y|<1\}$, as before, and let $u$ be the harmonic extension of the boundary data $u(x,\pm1)=|x|$, with at most linear growth. The Poisson kernel at the center line is
\[
P_0(t)=\frac12\operatorname{sech}\!\left(\frac{\pi t}{2}\right),
\qquad \int_\R P_0(t)\,dt=1.
\]
More generally, the density of the horizontal coordinate of planar Brownian motion exiting $\Pi$, starting at $(x,y)$, is $P_y(t-x)$, where
\[
\begin{aligned}
P_y(t)&=\frac14\left[
\frac{\cos(\pi y/2)}{\cosh(\pi t/2)-\sin(\pi y/2)}+
\frac{\cos(\pi y/2)}{\cosh(\pi t/2)+\sin(\pi y/2)}\right]\\
&=\frac{\cos(\pi y/2)\cosh(\pi t/2)}{2(\cosh^2(\pi t/2)-\sin^2(\pi y/2))}.
\end{aligned}
\]
In particular,
\begin{equation*}
u(x,y)=\int_\R |x+t|P_y(t)\,dt,
\qquad |y|<1.
\end{equation*}
(The kernel above is the sum of the two boundary hitting densities.)

An equivalent analytic description is especially useful. Put
\begin{equation}
\label{eq:F}
F(z)=\frac4\pi\arctan\left(\tanh\frac{\pi z}{4}\right),\qquad z\in \Pi,
\end{equation}
with the analytic branch satisfying $F(0)=0$. Then
\[
F'(z)=\operatorname{sech}\left(\frac{\pi z}{2}\right).
\]
For $x\ne0$, the boundary values satisfy
\[
\Re F(x\pm i)=\operatorname{sgn}(x).
\]
Indeed, the boundary derivative is the tangential derivative of $|x|$. By symmetry and uniqueness of the harmonic extension (in the linear-growth class),
\begin{equation}\label{eq:gradient}
u_x(x,y)-iu_y(x,y)=F(x+iy).
\end{equation}
One explicit formula for $u$ is therefore
\begin{equation*}
u(x,y)=a+\Re\int_0^{x+iy}F(\zeta)\,d\zeta,
\end{equation*}
where the integral is path independent in $\Pi$ and $a=u(0,0)$.

Recall from \eqref{eq:a} that $a=u(0,0)=8\Gcat/\pi^2$.

Define, for $0\le s<1$,
\[
q(s):=\frac4\pi\operatorname{arctanh}\left(\tan\frac{\pi s}{4}\right),
\qquad
\Phi(t):=\int_0^{\min(t,1)}q(s)\,ds\quad(t\ge0).
\]
For $0<y<1$, \eqref{eq:gradient} and $\arctan(iw)=i\operatorname{arctanh}w$ give
\[
F(iy)=iq(y),\qquad u_y(0,y)=-q(y).
\]
The second equality follows from \eqref{eq:gradient}  and evenness in $x$, which implies $u_x(0,y)=0$.
Hence
\begin{equation}
\label{eq:axis}
u(0,y)=a-\Phi(|y|),\qquad |y|\le1.
\end{equation}
Since $u(0,1)=0$, it follows that
\begin{equation}
\label{eq:phi-one}
\Phi(1)=a=\frac{8\Gcat}{\pi^2}.
\end{equation}

We now use $\Phi$ to construct a superharmonic majorant for the obstacle
\[ 
O(x, y):= \Phi(|y|) -|x|\,.
\]
Define on $\R^2$
\begin{equation*}
U(x,y)=\begin{cases}
a-u(x,y),& |y|<1,\\
a-|x|,& |y|\ge1
\end{cases}
= a \mathcal G(x, y).
\end{equation*}
We verify the properties
\begin{equation*}
\Delta U\le0\quad\text{in distributions},\qquad
U(x,y)\ge\Phi(|y|)-|x|,\qquad U(x,0)\le0.
\end{equation*}

The estimate proved in Subsection~\ref{Bell1}, together with
\eqref{eq:axis}, gives
\[
u(x,y)-|x|\le u(0,y)=a-\Phi(|y|),\qquad |y|<1.
\]
Thus $U(x,y)\ge\Phi(|y|)-|x|$ inside $\Pi$.
Outside $\Pi$, equality holds because $\Phi(|y|)=a$.

Since $U=a\mathcal G$ and $a>0$, the superharmonicity established
in Subsection~\ref{Bell1} gives $\Delta U\le0$ in distributions.
Moreover, \eqref{H00} gives
\begin{equation}
\label{at0}
U(x,0)=a\mathcal G(x,0)\le0.
\end{equation}

\subsection{The integral estimate for harmonic conjugates}
\label{hv-one-dimensional-section}
\begin{theorem}
Let $g\in L^1(\T)$ be real valued, and let $h$ be the almost-everywhere boundary value of its harmonic conjugate, normalized by $h(0)=0$ in the disk. For every $\lambda>0$,
\begin{equation}\label{eq:main}
\int_\T\Phi\left(\frac{|h(e^{it})|}{\lambda}\right)\,dm(t)
\le\frac1\lambda\int_\T |g(e^{it})|\,dm(t).
\end{equation}
\end{theorem}
\begin{proof}
First suppose $g$ is a real trigonometric polynomial or a smooth function on $\T$. Let $f=g+ih = g+i\cH g$ be analytic in the disk,  with $h(0)=0$. Since $U$ is superharmonic and analytic maps are conformal wherever their derivatives do not vanish, $z\mapsto U(g(z)/\lambda,h(z)/\lambda)$ is superharmonic. The mean-value inequality gives (measure on $\T$ is normalized)
\[
\int_\T
U\left(\frac{g(re^{it})}{\lambda},\frac{h(re^{it})}{\lambda}\right)\,dm(t)
\le U\left(\frac{g(0)}\lambda,0\right)\le 0.
\]
The last inequality is just \eqref{at0}.

By the obstacle domination property,
\[
\int_\T\Phi\left(\frac{|h(re^{it})|}{\lambda}\right)\,dm(t)
\le\frac1\lambda\int_\T|g(re^{it})|\,dm(t).
\]
Let $r\uparrow1$. For general $g\in L^1$, apply this to its Poisson regularizations (or Fej\'er approximants), whose conjugates converge almost everywhere along a suitable sequence. Fatou's lemma and $L^1$ convergence of the real parts give \eqref{eq:main}.
\end{proof}

Since $\Phi$ is increasing and $\Phi(t)=a$ for $t\ge1$, \eqref{eq:main} implies
\[
a\,|\{t:|h(e^{it})|\ge\lambda\}|
\le\frac{\|g\|_{L^1(\T)}}\lambda.
\]
Using \eqref{eq:phi-one}, we recover the weak-type estimate:
\begin{equation*}
|\{t:|h(e^{it})|\ge\lambda\}|
\le\frac{\pi^2}{8\Gcat}\frac{\|g\|_{L^1(\T)}}\lambda.
\end{equation*}

The $\Phi$ estimate \eqref{eq:main} gives more than the sharp weak-type estimate. Since $\Phi(0)=\Phi'(0)=0$ and $\Phi''(t)=\sec(\pi t/2)\ge1$ for $0\le t<1$, we see that on $t\in [0,1]$
\[
\Phi(t) \ge t^2/2\,.
\]
Thus \eqref{eq:main} also gives
\begin{equation*}
a \lambda |\{ |\cH g| >\lambda\}|+\frac1{2\lambda} \int_{|\cH g| <\lambda} |\cH g|^2 \le \|g\|_1\,.
\end{equation*}

\subsection{The integral estimate for horizontal Riesz transforms}
\label{sec:horizontal-integral}

With $u$ extended by $u(x,y)=|x|$ outside the strip, put
\begin{equation*}
 a=u(0,0)=\frac{8\Gcat}{\pi^2},\qquad
 U=a-u,\qquad \Bell(x,y)=U(x,y)+|x|.
\end{equation*}
The function $\Phi$ satisfies
\begin{equation*}
 \Phi(t)=a-u(0,\min\{t,1\}),\qquad t\ge0.
\end{equation*}
The preceding arguments give
\[
 0\le\Phi(|y|)\le\Bell(x,y)\le a,\qquad
 \Phi(t)=a\ (t\ge1),\qquad
 \Phi(t)\ge t^2/2\ (0\le t\le1).
\]
Here the upper bound $\Bell\le a$ follows from $u\ge|x|$.
The proof of Theorem~\ref{thm:horizontal} will establish the
following stronger estimate.

\begin{theorem}\label{thm:horizontal-integral}
Under the assumptions of Theorem~\ref{thm:horizontal}, for every
real $f\in L^1(X)$ and every $\lambda>0$,
\begin{equation*}
 \lambda\int_X
 \Bell\!\left(\frac{f}{\lambda},\frac{Rf}{\lambda}\right)\,d\mu
 \le\|f\|_1.
\end{equation*}
Consequently,
\begin{equation}\label{eq:main-distribution}
 \lambda\int_X\Phi\!\left(\frac{|Rf|}{\lambda}\right)\,d\mu
 \le\|f\|_1,
\end{equation}
and
\begin{equation}\label{eq:main-truncated}
 \cD\lambda\,\mu\{|Rf|\ge\lambda\}
 +\frac1{2\lambda}\int_{\{|Rf|<\lambda\}}|Rf|^2\,d\mu
 \le\|f\|_1.
\end{equation}
\end{theorem}

The pointwise bounds for $\Phi$ give
\eqref{eq:main-truncated}. The proof of
Theorem~\ref{thm:horizontal-integral} is completed in
Proposition~\ref{tr-upper}.

\section{Proofs of the main theorems}\label{sec:finite}

We first prove the Bellman estimate for a finite transverse space.

\begin{theorem}\label{fp-finite}
Let $E$ be finite with positive weights of total mass one, and let
$P_t=e^{-tL}$ be symmetric and sub-Markov. For every real
$f\in L^2(\T\times E)$,
\begin{equation}\label{fp-finite-U}
\int_{\T\times E}U(f,Rf)\,d\mu\le0.
\end{equation}
Equivalently, for every $\lambda>0$,
\begin{equation}\label{fp-finite-Bellman}
\lambda\int_{\T\times E}
\Bell\!\left(\frac f\lambda,\frac{Rf}\lambda\right)\,d\mu
\le\|f\|_1.
\end{equation}
The constant is independent of $E$, its weights and $L$.
\end{theorem}

We prove the estimate under a strip restriction, then remove that
restriction by regularization, maximization, and confinement.

\subsection{The estimate under a strip restriction}
\label{ndim}\label{fp-finite-section}

Recall the function $F$ from \eqref{eq:F}. For $w\in \Pi$, choose the branches in
\begin{equation*}
 F(w)=\frac4\pi\arctan\!\left(\tanh\frac{\pi w}{4}\right),
 \qquad F(0)=0.
\end{equation*}
Then $ |\operatorname{Re}F|<1$ and
\begin{equation}
\label{sech}
 \qquad
 F'(w)=\operatorname{sech}(\pi w/2),
 \qquad
 \operatorname{Re}F'(w)>0.
\end{equation}
The function $F$ is the composition $-iS_*(iz)$ of $S_*(z)$, $z\in \bD$, having boundary values 
\[
S_*(e^{ix}) = \frac2\pi\log\left|\cot\frac x2\right| +i \operatorname{sgn}(\sin x)
\]
and the conformal map $\Pi\to\bD$
\[
z= \tanh\frac{\pi w}{4}\,,\qquad F(w)=-iS_*\!\left(i\tanh\frac{\pi w}{4}\right).
\]
It is also 
\[
F= 2\pd_z u
\]
where $ u$ was constructed in the Bellman function proof of Davis' inequality above.

\medskip

Let $\Psi'=F$, $\Psi(0)=0$.

The function $\Psi$ is $u+i\tilde u$ up to an additive constant. 
Since
\[
 \partial_u\{u\operatorname{Re}F(u+iv)- \operatorname{Re}\Psi(u+iv)\}
 =u\operatorname{Re}F'(u+iv),
\]
$u\operatorname{Re}F(u+iv)-\operatorname{Re}\Psi(u+iv)$ is minimized at $u=0$.  The identity
\[
 F(is)=\frac{4i}{\pi}
 \operatorname{arctanh}\!\left(\tan\frac{\pi s}{4}\right)
\]
gives, after integration along the imaginary axis with respect to $d(iv)$,
\begin{equation}
\label{eq:strip}
 u\operatorname{Re}F(u+iv)-\operatorname{Re}\Psi(u+iv)\geq\Phi(|v|),
 \qquad |v|<1,
\end{equation}

We extend $F$ vertically beyond the strip.
For $|v|<1$ or $u\ne0$ define
\[
\widetilde F(u+iv)=
\begin{cases}
F(u+iv),& |v|<1,\\
\displaystyle\lim_{t\uparrow1}F(u+it),&v\ge1,\\
\displaystyle\lim_{t\downarrow-1}F(u+it),&v\le-1.
\end{cases}
\]

\[
\operatorname{Re}\widetilde F(u+iv)=\operatorname{sgn}(u), \quad |v|\ge 1\,.
\]

At $u=0$, $|v|\ge1$, set $\Re\widetilde F(u+iv)=0$; only its real part is used there.

We define the real-valued extension by vertical freezing:
\begin{equation}
\label{extPsi}
\operatorname{Re}\widetilde\Psi(u+iv)
=\begin{cases}
\operatorname{Re}\Psi(u+iv),&|v|<1,\\
|u|-\Phi(1),&|v|\ge1.
\end{cases}\quad
\end{equation}

Let
\[
\widetilde H(u,v)=u\operatorname{Re}\widetilde F(u+iv)
- \operatorname{Re}\widetilde\Psi(u+iv)
\]
where $\widetilde F$ is defined, and assign its continuous value
$\Phi(1)$ on $u=0$, $|v|\ge1$. Then
\[
\widetilde H(u,v)=\Phi(1)=\Phi(|v|),
\qquad |v|\ge1.
\]
Together with \eqref{eq:strip}, this gives the inequality everywhere:
\begin{equation}
\label{main}
 u\operatorname{Re}\widetilde F(u+iv)-\operatorname{Re}\widetilde\Psi(u+iv)\geq\Phi(|v|).
\end{equation}

Recall the extension of $\mathcal G_C$ outside the strip $\Pi$ in \eqref{eq:paste}: $u$ was extended by $|x|$, and the resulting $\mathcal G_C$ is superharmonic. Thus the extension of $u$ is subharmonic. The same holds for the real-valued extension of $\Psi$
in \eqref{extPsi}, since $\operatorname{Re}\Psi$ and $u$ differ by a constant.

This subharmonicity means in particular 
\begin{equation*}
\int_\T \operatorname{Re}\widetilde\Psi(h(\zeta) +i g(\zeta))\,dm(\zeta) \ge 0
\end{equation*}
if $g$ is a smooth function on $\T$ and $h=-\cH g$ and $\int_\T g\,dm=0$. We used here that $\operatorname{Re}\widetilde\Psi(0) =0=\operatorname{Re}\Psi(0)$ and the subharmonicity established above.

Substitute $u= h(\zeta)$, $v=g(\zeta)$ into \eqref{main} and integrate over $\T$ to obtain
\begin{equation}
\label{Phi2}
\int_\T \Phi(|g|) \le \int_\T h\,  \Re \widetilde F (h+ig)
\end{equation}
The same inequality holds for $H^2$ boundary data by radial approximation: $h\Re\widetilde F(h+ig)$ is continuous in $(h,g)$ and bounded in absolute value by $|h|$.

We first treat the conservative case. Suppose now that $g$ also depends on a parameter $a\in E$. On $E$, which we take to be finite for the moment, choose weights
$\nu_a>0$ with
$\sum_a\nu_a=1$, and let $-L$ be the generator of a reversible
continuous-time Markov chain on $E$. Thus
\[
 P_t\alpha(a)=\sum_{b\in E}p_t(a,b)\alpha(b)\nu_b,
 \qquad
 p_t(a,b)=p_t(b,a),
 \qquad
 \sum_b p_t(a,b)\nu_b=1.
\]
All integrals on $\mathbb T\times E$ use normalized Haar measure in
$\T$ and the probability measure $\nu$ in the second variable. 

Let $L$ be a Laplacian corresponding to this chain.
Let $P_t=e^{-tL}$ be a strongly continuous symmetric Markov
semigroup on $L^2(E,\nu)$. Thus $L$ is nonnegative and self-adjoint,
$P_t$ preserves positivity, and its consistent $L^\infty$ realization
satisfies $P_t1=1$. On either
$\mathbb X=\mathbb T\times E$ or $\mathbb R\times E$, put
\[
 A=(-\partial_x^2+L)^{1/2},\qquad R=-\partial_xA^{-1}.
\]
We use normalized Haar measure $dx/(2\pi)$ on
$\mathbb T=\mathbb R/(2\pi\mathbb Z)$ and Lebesgue measure on
$\mathbb R$; $\mu$ denotes the resulting product measure with $\nu$.
We first take $E$ finite. Bounds independent of $E$ allow passage to
general transverse spaces, as proved in Subsection~\ref{tr-section}.

\medskip

 The operator $R$ is the bounded joint spectral multiplier on
$L^2$, with value zero on $\ker A$. On the circle,
\begin{equation*}
 \widehat{Rf}(m)=-im(m^2+L)^{-1/2}\widehat f(m)
 \quad(m\ne0),\qquad \widehat{Rf}(0)=0.
\end{equation*}
On the line replace $m$ by $\xi$ and use
$\widehat f(\xi)=\int e^{-ix\xi}f(x)\,dx$.
In particular, $\|R\|_{2\to2}\leq1$ and
$R^*=-R$. We write $R_L$ when the transverse generator must be specified.

The generator $A^2=-\partial_x^2+L$ has the product semigroup
\[
 e^{-tA^2}=e^{t\partial_x^2}\otimes P_t.
\]
Thus $R$ is the coordinate Riesz transform of this product generator
in the ordinary line or circle direction. We call this direction
\emph{horizontal} and the $E$ variables \emph{transverse}. The
first-order operator is supplied by the line or circle factor, so
the construction requires no derivatives in the $E$ variables.

The question here is whether adjoining independent Markov dynamics
preserves the exact Hilbert-transform weak endpoint constant.
The square root couples horizontal frequency to the transverse
variables, yet the exact norm remains the same.
We prove this first for the finite Markov chain $L$.

\medskip

Let $f\in L^2(\mathbb T\times E;\mathbb R)$. We impose the strip
restriction below, when proving the weighted positivity estimate.  
Put
\[
 g=Rf,
 \qquad
 h=-\cH g.
\]
For each $a\in E$, the function $h(\cdot,a)+ig(\cdot,a)$ is the
boundary value of an $H^2$ function $S_a$ with $S_a(0)=0$.  
Put
\[
 q(\cdot,a)=\operatorname{Re}\widetilde F(S_a).
\]
By the definition of $\operatorname{Re}\widetilde F$,
\begin{equation*}
|q|\le 1
\end{equation*}
For every $a\in E$, \eqref{Phi2} gives
\[
\int_\T \Phi(|g(\cdot, a)|) \le \int_\T q(\cdot, a) h(\cdot, a)\,.
\]
Summing over the states,
\begin{equation*}
\int_{\T\times E} \Phi(|g|) \le \langle q, h\rangle\,.
\end{equation*}
 
We next prove that $\langle q, Th\rangle\ge 0$ for an operator $T$ related to the horizontal Riesz transform.

For $t>0$, let $H_t$ be the horizontal multiplier
\[
 \widehat{H_tu}(m)=|m|^{-1}e^{-tm^2}\widehat u(m),
 \qquad m\neq0,
\]
with value zero at $m=0$.  Its multipliers decrease with $|m|$.
Reversibility and the Markov identity give the equality below; the inequality will be proved under the strip assumption \eqref{ImS}:
\begin{align}
  \label{eq:Markov-positive}
 \langle q,H_t(I-P_t)h\rangle
 &=\frac12\sum_{a,b}\nu_a\nu_bp_t(a,b)
 \langle q_a-q_b,H_t(h_a-h_b)\rangle_{L^2(\mathbb T)}
 \geq 0.                                         
\end{align}

If we prove
\begin{equation}
\label{eq:weighted-positive}
 \langle q_a-q_b,H_t(h_a-h_b)\rangle_{L^2(\mathbb T)}\ge 0\,.
 \end{equation}
 then \eqref{eq:Markov-positive} is proved.
 
Let $|D_x|$ denote the horizontal Fourier multiplier $|m|$.
On the subspace of functions with zero horizontal mean, set
\[
 T=\frac{(-\partial_x^2+L)^{1/2}}{|D_x|}-I.
\]
For $r>0$ and $s\geq0$ the following identity holds
\[
 \sqrt{1+\frac{s}{r^2}}-1
 =\frac{1}{2\sqrt\pi\,r}\int_0^\infty
 e^{-tr^2}(1-e^{-ts})t^{-3/2}\,dt.
\]
The operator $T$ is bounded because $|m|\geq1$ and
$L$ is a finite matrix. By the preceding identity, the following integral converges in operator norm and represents $T$:
\begin{equation}
\label{eq:T-representation}
 T=\frac1{2\sqrt\pi}\int_0^\infty
 H_t(I-P_t)t^{-3/2}\,dt.
\end{equation}
Indeed $\|I-P_t\|_{2\to2}\leq\min(t\|L\|,2)$ and
$\|H_t\|_{2\to2}\leq e^{-t}$.

Under the strip assumption \eqref{ImS}, we see that the previous formula and \eqref{eq:Markov-positive} imply
\begin{equation}
\label{eq:qTh}
 \langle q,T h\rangle\geq0.
\end{equation}

\medskip

It remains to prove \eqref{eq:weighted-positive}.

Let
$S_1=h(\cdot,a)+ig(\cdot,a), S_2=h(\cdot,b)+ig(\cdot,b)\in H^2$
have zero value at the origin, and let $h_j$
be the real boundary parts of $S_j$ and $q_j$ be the boundary parts of
$\Re\widetilde F(S_j)$.
We now assume that
\begin{equation}
\label{ImS}
|\Im S_1|\le1,\quad |\Im S_2|\le1
\quad\Longleftrightarrow\quad
|Rf(\cdot,a)|\le1,\quad |Rf(\cdot,b)|\le1
\end{equation}
on the boundary. The maximum principle puts $S_1,S_2$ in the open
strip inside $\bD$, so the extension is not needed in the following
holomorphic calculation.
If $M_w$ is the multiplier $w_{|m|}$ on nonzero frequencies such that
$w_1\geq w_2\geq\cdots\geq0$, then we wish to see that
\begin{equation}
\label{eq:weighted-positive1}
\langle q_1-q_2,M_w(h_1-h_2)\rangle_{L^2(\mathbb T)}\geq0.
\end{equation}
Indeed, let $Q_j=F(S_j)$, $Q=Q_1-Q_2$, and $S=S_1-S_2$.
Then
\begin{equation*}
Q(z)=B(z)S(z),\qquad
B(z)=\int_0^1F'(S_2(z)+t(S_1(z)-S_2(z)))\,dt.
\end{equation*}
The function $B$ is holomorphic and
\begin{equation*}
\Re B\geq0
\end{equation*}
by \eqref{sech} and the convexity of the strip.

We have $h_1-h_2=(S+\overline S)/2$ and
$q_1-q_2=(Q+\overline Q)/2$. Since $S(0)=Q(0)=0$,
\[
\frac12\Re\sum_{j=1}^\infty\widehat Q(j)\overline{\widehat S(j)}
=\frac12\Re\int_\T Q\overline S
=\int_\T(q_1-q_2)(h_1-h_2).
\]
Here $Q_j\in H^2$, since $\Re F(S_j)$ is bounded and $F(S_j(0))=0$.
Thus the coefficient products are absolutely summable.

The identity $Q=B(S_1-S_2)$, together with
analyticity of $B$, gives for any $n$
\[
A_n:=\frac12\Re\sum_{j=1}^n
\widehat Q(j)\overline{\widehat S(j)}
=\lim_{\rho\uparrow1}\frac12\int_\T
\Re B(\rho\zeta)\,
|P_n(S_1-S_2)(\rho\zeta)|^2\,dm(\zeta)\geq0,
\]
where $P_n$ is the projection on $\operatorname{span}\{e^{ikx}\}_{k=1}^n$
and the integral is normalized. Indeed, the coefficients through
order $n$ of $BP_nS$ agree with those of $BS=Q$.
The radial limit retains any singular boundary mass of $\Re B$.

It follows that
\[
\begin{aligned}
\langle q_1-q_2,M_w(h_1-h_2)\rangle_{L^2(\mathbb T)}
&=\frac12\Re\sum_{j=1}^\infty
w_j\widehat Q(j)\overline{\widehat S(j)}\\
&=w_1A_1+w_2(A_2-A_1)+w_3(A_3-A_2)+\cdots\\
&=w_\infty A_\infty+
\sum_{n=1}^\infty(w_n-w_{n+1})A_n\geq0,
\end{aligned}
\]
where $w_\infty=\lim_nw_n$ and $A_\infty=\lim_nA_n$.

This proves \eqref{eq:Markov-positive} under the
assumption \eqref{ImS}. Under this assumption in every state we have
\begin{equation*}
\langle q,h\rangle\geq\int_{\T\times E}\Phi(|Rf|),
\qquad \langle q,Th\rangle\geq0,
\end{equation*}
where the pairings are taken on $\T\times E$.

But $h+Th = (-\pd_x^2 + L)^{1/2}|D_x|^{-1}[-\cH Rf]= f-\bar f$, where $\bar f$ is the horizontal average in each state. Under the strip assumption, analyticity gives $\bar q=0$. Hence $\langle q,h+Th\rangle=\langle q,f\rangle\le\|f\|_1$, since $|q|\le1$.
This proves the following estimate.
\begin{lemma}
\label{Rfsmall}
If for every $a\in E$ we have $|Rf(\cdot, a)|\le 1$ then
\begin{equation*}
\int_{\T\times E} \Phi(|Rf|) \le \|f\|_1\,.
\end{equation*}
\end{lemma}
To remove the strip restriction, we shall use the following
stronger estimate:
\begin{equation}
\label{less1}
\int U(f, Rf)\le 0\,\, \text{if}\,\, |Rf|\le 1\,\, \text{everywhere}\,.
\end{equation}

\medskip

Suppose $g=Rf$ and $\|g\|_\infty<1$. Put $h=-\cH g$ and
$q=\Re F(h+ig)$, state by state. A bar denotes horizontal averaging, and $D=|D_x|$.
Analyticity and radial limits give $\bar q_a=0$ and
$\int_{\mathbb T}U(h_a,g_a)\,dm=0$.
Concavity in the first coordinate yields
\[
U(f,g)\le U(h,g)-q(f-h),\qquad
f-h=(AD^{-1}-I)h+\bar f.
\]
On functions with zero horizontal mean, $AD^{-1}-I=T$.
The representation \eqref{eq:T-representation} and the positivity
\eqref{eq:qTh} give $\langle q,Th\rangle\ge0$ under the strip restriction.
Using $\bar q=0$ and $\int U(h,g)\le 0$ in the concavity inequality, we obtain
\[
\int U(f,Rf)\le-\langle q,Th\rangle\le0.
\]
The case $|Rf|\le1$ follows by applying the strict-strip case to $rf$ and letting $r\uparrow1$, using continuity and $|U(rf,rRf)|\le a+|f|$.

\paragraph{Finite semigroups with killing.}
The same strip estimate holds when $P_t1\le1$. Indeed, put
\[
 d_t(a)=1-\sum_b p_t(a,b)\nu_b\ge0.
\]
Then the reversible pairing used above becomes
\[
\begin{split}
\langle q,H_t(I-P_t)h\rangle
={}&\frac12\sum_{a,b}\nu_a\nu_b p_t(a,b)
 \langle q_a-q_b,H_t(h_a-h_b)\rangle\\
&+\sum_a\nu_a d_t(a)\langle q_a,H_t h_a\rangle\ge0.
\end{split}
\]
The additional terms are nonnegative by
\eqref{eq:weighted-positive1}, comparing $S_a$ with zero.
The variational argument below uses only the nonnegativity and self-adjointness
of $L$, its boundedness on the finite space, and the sub-Markov
property of $e^{-sA}$. We therefore permit a symmetric sub-Markov transverse semigroup
on the finite weighted space throughout the rest of the finite-space proof.

\subsection{Regularization and maximization}\label{fp-remove-section}

We remove the restriction $|Rf|\le1$ by a variational argument. 

We consider the extremal problem
\[
\int_{\T\times E} U(\phi, R\phi)\,d\mu\to \text{maximize}\,.
\]
The idea is to obtain a maximizer $f:= \phi^{max}$ by upper semicontinuity and prove that 
\[
|Rf|\le 1\,.
\]
Then the strip estimate \eqref{less1} gives
\[
\int_{\T\times E} U(f,Rf)\,d\mu\le0.
\]
By maximality, the same inequality would hold for every input.
We make this argument rigorous through regularization. 
We use the following properties of $U$:
\begin{align*}
& U(x, y) \,\, \text{ is even and concave in} \,\, x\,\, \text{globally};
\\
& U(x, y) \,\, \text{has uniform in }\,\, x\,\,\text{modulus of continuity in }y\,\, \omega;
\\
& U(x, y)= a-|x|, \,\, |y|\ge 1.
\end{align*}

Fix $\delta,\eps>0$. We replace $U$ by $U_\delta$ such that
\begin{align*}
& \|U_\delta -U\|_\infty \le \omega(\delta);
\\
& U_\delta\,\, \text{ is globally superharmonic};
\\
&U_\delta (x, y) = a-|x|,\,\, |y|\ge 1+\delta.
\end{align*}
For this, convolve $U(x, y)$ vertically with an even nonnegative smooth probability kernel in $y$ supported on $[-\delta, \delta]$.

\medskip

Consider the regularized variational problem
\[
J_{\delta, \eps} (\phi):=\int_{\T\times E} U_\delta(\phi, R\phi)\,d\mu-\eps\int_{\T\times E} |R\phi|^2\,d\mu\to \text{maximize}\,.
\]

The maximization is over real $L^2(\T\times E)$.

We prove the properties used in the variational argument, including the
uniform continuity asserted above. The arguments below also
apply when the finite transverse semigroup is sub-Markov.

\begin{lemma}\label{fp-strip-properties}
The function $U$ is continuous, superharmonic, even in each variable,
concave and one-Lipschitz in $x$. It satisfies
$-|x|\le U(x,y)\le a-|x|$. The function
$\Bell(x,y)=U(x,y)+|x|$ is continuous and
\begin{equation*}
0\le\Phi(|y|)\le\Bell(x,y)\le a.
\end{equation*}
Moreover,
\begin{equation*}
 \omega(\delta):=\sup_{\substack{x,y\in\R\\|h|\le\delta}}
 |U(x,y+h)-U(x,y)|\longrightarrow0\qquad(\delta\downarrow0).
\end{equation*}
\end{lemma}

\begin{proof}
The first properties follow from the Poisson formula and the preceding
calculations. Put $e(x,y)=u(x,y)-|x|$, with $e=0$ outside $\Pi$.
Convexity in $x$, harmonicity and evenness in $y$ give
\[
 0\le e(x,y)\le e(x,0)
 =2\int_{|x|}^\infty(t-|x|)P_0(t)\,dt\longrightarrow0
 \qquad(|x|\to\infty).
\]
Also $e(x,y)\le u(0,y)\to0$ as $|y|\uparrow1$, uniformly in $x$.
Thus $e\in C_0(\R^2)$ is uniformly continuous, and
$U(x,y)=a-|x|-e(x,y)$ gives the assertion.
\end{proof}

Let $\zeta_\delta$ be the even nonnegative smooth probability kernel used
above, so that
\[
U_\delta(r,s)=\int_\R U(r,s-t)\zeta_\delta(t)\,dt,
\qquad \operatorname{supp}\zeta_\delta\subset[-\delta,\delta].
\]

\begin{lemma}
\label{fp-regularization}
The function $U_\delta$ is superharmonic, even and concave in $r$,
and
\begin{equation*}
 \|U_\delta-U\|_\infty\le\omega(\delta),\qquad
 -|r|\le U_\delta(r,s)\le a-|r|,\qquad
 U_\delta(r,s)=a-|r|\quad(|s|\ge1+\delta).
\end{equation*}
It is one-Lipschitz in $r$.  The function $\partial_sU_\delta$ is bounded and
globally Lipschitz in $(r,s)$, and $\partial_sU_\delta(r,s)=0$ for
$|s|\ge1+\delta$.  Also $U_\delta(0,0)\ge0$.
\end{lemma}

\begin{proof}
Convolution preserves symmetry, concavity and superharmonicity.
Lemma~\ref{fp-strip-properties} gives the bounds and exterior formula.
When differentiating in $s$, use $U(r,s)+|r|\in[0,a]$. Thus
\[
 \|\partial_sU_\delta\|_\infty\le a\|\zeta_\delta'\|_1,\qquad
 \|\partial_s^2U_\delta\|_\infty
       \le a\|\zeta_\delta''\|_1.
\]
The one-Lipschitz bound for $U$ in $r$ also gives
$|\partial_sU_\delta(r,s)-\partial_sU_\delta(r',s)|
\le\|\zeta_\delta'\|_1|r-r'|$.
Finally $U(0,s)=\Phi(|s|)\ge0$, so $U_\delta(0,0)\ge0$.
\end{proof}

Let $S_t=e^{-tD}$ be horizontal Poisson smoothing, with kernel
$\mathsf P_t$ relative to normalized Haar measure $m$, and set $K=R-\cH$.

\paragraph{Poisson comparison.}

For $R=\cH$, superharmonicity gives Poisson comparison directly.
The multiplier of $K$ decays quadratically, giving an error of order
$t$. The gain from the quadratic penalty will supply the half
derivative needed for the contact test.

\begin{lemma}
\label{fp-Poisson-comparison}
There is a constant $C_{\delta,L}$ such that for every real
$f\in L^2(\T\times E)$ and $t>0$,
\begin{equation}\label{fp-Poisson-estimate}
 J_{\delta,\eps}(S_tf)\ge J_{\delta,\eps}(f)
 +\eps\bigl(\|Rf\|_2^2-\|S_tRf\|_2^2\bigr)
 -C_{\delta,L}t\bigl(\|f\|_2+\|f\|_2^2\bigr).
\end{equation}
The constant is independent of $\eps$, $t$, and $f$.
\end{lemma}

\begin{proof}
At horizontal frequency $m\ne0$, the spectral theorem gives
\begin{equation*}
 \|\widehat K(m)\|_{L^2(E)\to L^2(E)}
 \le\frac{\|L\|}{2m^2}.
\end{equation*}
This follows from
$0\le1-(1+\sigma/m^2)^{-1/2}\le \sigma/(2m^2)$ for $\sigma\ge0$.
The zero-frequency multiplier is zero.
In particular, for $k=Kf$,
$\|D k\|_2+\|D^{1/2}k\|_2\le C_L\|f\|_2$.

For each state and each fixed $x$, apply superharmonic Poisson
comparison to the analytic extension of $f+i\cH f$, keeping the
imaginary translation $S_tk(x)$ fixed.  This gives
\begin{equation}\label{fp-pointwise-Poisson}
 U_\delta(S_tf(x),S_tRf(x))
 \ge\int_\T\mathsf P_t(x-y)
 U_\delta\bigl(f(y),\cH f(y)+S_tk(x)\bigr)\,dm(y).
\end{equation}
To justify this at the cusp, first convolve $U_\delta$ in $r$.
Apply superharmonic comparison on smaller disks and pass to radial
boundary values, using the bound by a constant plus $|f|$ and radial
$L^2$ convergence. The one-Lipschitz estimate in $r$ removes the
extra convolution uniformly.

Put $\eta=\partial_sU_\delta(f,Rf)$ and expand the second coordinate of the
integrand in \eqref{fp-pointwise-Poisson} about $Rf(y)$.
After integration in $x$, the linear term is
$\langle\eta,(S_{2t}-I)k\rangle$.
For each state the quadratic error is bounded using the exact identity
\begin{equation*}
 \iint_{\T^2}\mathsf P_t(x-y)|S_tk(x)-k(y)|^2\,dm(x)\,dm(y)
 =\|k\|_2^2-\|S_tk\|_2^2.
\end{equation*}
The cross term is $2\langle S_tk,S_tk\rangle$.
After summing states, the Taylor remainder is bounded by
\[
 \frac12\|\partial_s^2U_\delta\|_\infty
          \bigl(\|k\|_2^2-\|S_tk\|_2^2\bigr).
 \]
The multiplier inequalities $1-e^{-2t|m|}\le2t|m|$ give
\begin{align*}
 |\langle\eta,(S_{2t}-I)k\rangle|
 &\le2t\|\eta\|_2\|D k\|_2
 \le C_{\delta,L}t\|f\|_2,\\
 \|k\|_2^2-\|S_tk\|_2^2
 &\le2t\|D^{1/2}k\|_2^2
 \le C_Lt\|f\|_2^2.
\end{align*}
Here $\|\eta\|_2\le\|\partial_sU_\delta\|_\infty$ because $\mu(\T\times E)=1$.
Add the quadratic penalties and use $RS_t=S_tR$ to obtain
\eqref{fp-Poisson-estimate}.
\end{proof}

\paragraph{Existence and regularity.}

Because $E$ is finite and $L$ is bounded, the spaces
\[
 W^{1/2,2}=W^{1/2,2}(\T;L^2(E,\nu))
 \quad\text{and}\quad \operatorname{Dom}(A^{1/2})
\]
coincide, with equivalent graph norms. We use this identification below.

\begin{proposition}
\label{fp-maximizer}
The functional $J_{\delta,\eps}$ attains its supremum on real $L^2(\T\times E)$.
Every maximizer $f$, with $g=Rf$, satisfies
$f,g\in W^{1/2,2}$.
\end{proposition}

\begin{proof}
Lemma~\ref{fp-regularization} gives
\[
 J_{\delta,\eps}(0)\ge0,\qquad
 J_{\delta,\eps}(f)\le a-\|f\|_1-\eps\|Rf\|_2^2.
\]
For the horizontal mean $\overline f$,
\begin{equation}\label{fp-inverse-bounds}
 \|f-\overline f\|_2\le\sqrt{1+\|L\|}\,\|Rf\|_2,\qquad
 \|\overline f\|_2\le(\min_j\nu_j)^{-1/2}\|f\|_1.
\end{equation}
The first is the inverse multiplier bound on $m\ne0$; the second
compares the weighted $\ell^2$ and $\ell^1$ norms of the horizontal
averages. Thus maximizing sequences with values bounded below are
bounded in $L^2$.

Suppose $f_\ell\rightharpoonup f$ in $L^2$.
For fixed $t>0$, $S_t$ is compact because $E$ is finite, so
$S_tf_\ell\to S_tf$ strongly.  The Lipschitz bounds for $U_\delta$
and boundedness of $R$ show that $J_{\delta,\eps}$ is strongly continuous.
Dropping the nonnegative energy gain in
\eqref{fp-Poisson-estimate} therefore gives
\[
 \limsup_{\ell\to\infty}J_{\delta,\eps}(f_\ell)
 \le J_{\delta,\eps}(S_tf)+Ct,
\]
where $C$ is independent of $\ell$ and $t$.
Letting $t\downarrow0$ proves weak upper semicontinuity.
Weak compactness of a bounded maximizing sequence now gives a maximizer.

At such a maximizer, $J_{\delta,\eps}(S_tf)\le J_{\delta,\eps}(f)$.
Poisson comparison implies
\[
 \eps(\|g\|_2^2-\|S_tg\|_2^2)\le C_ft.
\]
Fourier expansion and Fatou's lemma yield
\begin{equation*}
 2\sum_{m\in\mathbb Z}|m|\|\widehat g(m)\|_{L^2(E)}^2
 \le\liminf_{t\downarrow0}
 \sum_m\frac{1-e^{-2t|m|}}t\|\widehat g(m)\|_{L^2(E)}^2
 \le\frac{C_f}{\eps}.
\end{equation*}
Hence $g\in W^{1/2,2}$.
The first multiplier bound in \eqref{fp-inverse-bounds}, now with
the factor $|m|$ in the Fourier sum, gives $f\in W^{1/2,2}$ as well.
Its horizontal mean is a vector in a finite-dimensional space and
causes no additional domain condition.
\end{proof}

Horizontal differentiation maps $W^{1/2,2}$ into $W^{-1/2,2}$.
We shall test the Euler identities with truncations and a Lipschitz
cutoff of $g$, using the following fact.

\begin{lemma}
\label{fp-composition}
Globally Lipschitz maps between finite-dimensional real spaces act
continuously by composition on $W^{1/2,2}(\T)$.
In particular, the same assertion holds statewise on $\T\times E$.
For $v\in W^{1/2,2}$, its truncations
$v_M=\max\{-M,\min\{v,M\}\}$ converge to $v$ in $W^{1/2,2}$.
\end{lemma}

\begin{proof}
The Lipschitz bound controls the Gagliardo difference quotient of
the composition by that of the input, and the finite measure controls
any nonzero value of the map at the origin.
To prove continuity, from any sequence converging in $W^{1/2,2}$
extract a subsequence whose $W^{1/2,2}$ errors are summable.
The sum of their absolute Gagliardo difference quotients belongs to
$L^2$ by Minkowski's inequality.
Adding the difference quotient of the limit gives an $L^2$ majorant
for the input differences along this subsequence, and hence for
the Lipschitz compositions.  Pointwise convergence and dominated
convergence prove strong convergence of their difference quotients
and of their $L^2$ values.  The subsequence criterion proves continuity
for the full sequence.
For truncations, use
$|v_M(x)-v_M(y)|\le|v(x)-v(y)|$ and pointwise convergence;
dominated convergence applies directly to the difference-quotient
and $L^2$ terms.
\end{proof}

\paragraph{Deduction from confinement.}

For the maximizer $f=f_{\delta,\eps}$ from Proposition~\ref{fp-maximizer}, we will prove
\begin{equation}
\label{Rfde}
|R f_{\delta, \eps}|\le 1+\delta\,.
\end{equation}

The confinement estimate will also give
\begin{equation*}
J_{\delta, \eps}(f_{\delta,\eps}) \le C(\delta, \eps)\to 0, \, \, \delta, \eps\to 0\,.
\end{equation*}
Assuming the confinement estimate, we obtain
\begin{align*}
&\int U_\delta (f_{\delta, \eps}, R f_{\delta, \eps})  \le \int U(f_{\delta, \eps}, R f_{\delta, \eps})  +\om(\delta) \le \notag
\\
& \int U\Big(\frac{f_{\delta, \eps}}{1+\delta}, \frac{R f_{\delta, \eps}}{1+\delta}\Big)  +2\om(\delta)  \le 2\om(\delta).
\end{align*}
The second inequality follows from evenness and concavity of $U$ in the first argument, its uniform modulus of continuity in the second, and the bound
\[
\Big| \frac{R f_{\delta, \eps}}{1+\delta} - Rf_{\delta, \eps}\Big| =\delta\frac{|Rf_{\delta, \eps}|}{1+\delta}\le \delta\,.
\]
where we used \eqref{Rfde}.

The third inequality follows from the strip estimate
\eqref{less1}, applied to $f_{\delta,\eps}/(1+\delta)$;
its transform has absolute value at most one by \eqref{Rfde}.

For any trigonometric polynomial $P(e^{ix}, a)$ on $\T\times E$, $J_{\delta, \eps}(P)$ is at most the maximum, so
\begin{align*}
&\int U_\delta(P, RP) \le  \int U_\delta (f_{\delta, \eps}, R f_{\delta, \eps})  - \eps \int |Rf_{\delta, \eps}|^2 + \eps \int |RP|^2\le
\\
& 2\om(\delta) + \eps \int |RP|^2\,.
\end{align*}
Let $\eps\to0$ and $\delta\to0$ to obtain
\begin{equation}
\label{UR}
\int \Phi(|Rf|) -\|f\|_1\le \int U(f, Rf) \le 0\,.
\end{equation}
The first inequality is the obstacle condition. The passage from trigonometric polynomials to real $L^2$ functions follows by $L^2$ approximation, boundedness of $R$, the one-Lipschitz dependence of $U$ on its first coordinate, and its uniform continuity in the second.

\bigskip

It remains to prove \eqref{Rfde}.

\subsection{The confinement estimate}
\label{fp-contact-section}

We suppress subscripts where possible: $J_{\delta,\eps}$ will be $J$. So let $f=f_{\delta,\eps}$ and put $g=Rf$.
We use variations of $J(f+h)-J(f)$ to obtain an Euler--Lagrange relation.
By Proposition~\ref{fp-maximizer}, a maximizer exists and
$f,g\in W^{1/2,2}=\operatorname{Dom}(A^{1/2})$.
Define the bad set
\[
B:=\{(x,a):|g(e^{ix},a)|>1+\delta\}.
\]
The measure $\mu$ remains the product measure on all of $\T\times E$; integration over the bad set is denoted by $\int_B\cdots\,d\mu$.

\paragraph{The Euler--Lagrange relation.}

For a convex functional \(\mathcal C\) on a Banach or Hilbert space \(X\), the \textbf{subdifferential}
\[
\partial \mathcal C(f)
\]
is the collection of continuous \textbf{linear} functionals \(\ell\in X^*\) such that
\begin{equation}
\label{ellh}
\mathcal C(f+h)\ge \mathcal C(f)+\ell(h)
\qquad\text{for every }h\in X.
\end{equation}
In a Hilbert space, one identifies \(\ell\) with an element \(p\in X\), writing
\[
\ell(h)=\langle p,h\rangle.
\]

The function \(U_\delta\) was smoothed only in the \(s\)-variable. In the exterior region,

\[
U_\delta(r,s)=a-|r|,\qquad |s|>1+\delta,
\]
so it still has a cusp at \(r=0\). Thus one cannot simply write \(\partial_rU_\delta(f,g)\).
Instead, freeze the second variable \(g\), perturb $f\to f+h$  and define the convex functional
\[
\mathcal C(f+h)=-\int U_\delta(f+h,g)\,d\mu .
\]
The Euler--Lagrange relation also involves variation in the second coordinate. Put
 \[
 \eta:=\partial_sU_\delta(f,g).
 \]
 Thus
\[
\begin{aligned}
J(f+h)-J(f)
={}&-\bigl(\mathcal C(f+h)-\mathcal C(f)\bigr)\\
&+\langle \eta-2\varepsilon g,Rh\rangle
+O(\|h\|_2^2).
\end{aligned}
\]
Since \(f\) maximizes \(J\), applying this with $th$, letting $t\downarrow0$, and using convexity gives
\begin{equation}
\label{subgrad}
R^*(\eta-2\varepsilon g)\in\partial\mathcal C(f).
\end{equation}

\begin{rem}
If $\mathcal C$ were linear $\langle c, \cdot\rangle$ we would write  Euler-Lagrange equation in the form
\[
R^*(\eta-2\varepsilon g) = c\,.
\]
\end{rem}

\medskip

But $\mathcal C$ is not linear, so using \(R^*=-R\) define
\begin{equation*}
p:=-R(\eta-2\varepsilon g).
\end{equation*}
Varying $h$ in \eqref{ellh}, the convex subgradient condition
\eqref{subgrad} gives pointwise:
\begin{equation}
\label{point}
p(r-f) \le U_\delta(f, g) - U_\delta(r,g)
\qquad\text{for every }r.
\end{equation}
Because \(U_\delta\) is one-Lipschitz in \(r\), every such subgradient of $-U_\delta(\cdot,g)$ has modulus at most \(1\). Therefore
\[
|p|\le1.
\]
We have two more Euler-Lagrange relations that follow from the definition of \(R\):
\[
Ap=\partial_x\eta-2\varepsilon\partial_xg,
\qquad
Ag=-\partial_xf.
\]
Indeed, since \(g=Rf=-\partial_xA^{-1}f\),
\[
p=-R(\eta-2\varepsilon g)
=\partial_xA^{-1}(\eta-2\varepsilon g).
\]

\paragraph{The localized test.}

Choose an even smooth function $\chi$ with $0\le\chi\le1$ and bounded derivative such that

\[
\chi(s)=0\quad\text{for }|s|\le 1+\delta,
\qquad
\chi(s)>0\quad\text{for }|s|>1+\delta.
\]
Then $\Theta(s)=\int_0^s\chi(r)\,dr$ is a smooth version of $
\operatorname{sgn}(s)(|s|-(1+\delta))_+.$
In particular,
\[
\Theta(s)=0\quad\Longleftrightarrow\quad |s|\le 1+\delta.
\]
We shall use the following two properties of $\Theta$:
\begin{equation}
\label{theta}
(s-r)(\Theta(s)-\Theta(r)) \ge (\Theta(s)-\Theta(r))^2,\quad r\Theta(r) \ge \Theta(r)^2
\end{equation}
The truncation
\[
f_M=\max(-M,\min(f,M))
\]
has the same sign as \(f\), and
\[
v_M:=f_M\chi(g)
\]
is essentially \(f\) localized to the bad region \(\{|g|>1+\delta\}\). The truncation is needed because \(f\) need not be bounded, while the map
$(r,s)\longmapsto r_M\chi(s)$ is globally Lipschitz and therefore acts continuously by composition on $W^{1/2,2}$.

On \(\{|g|>1+\delta\}\),
\[
U_\delta(r,g)=a-|r|.
\]

The concave superdifferential of \(a-|r|\) is

\[
\partial_r^+(a-|r|)
=
\begin{cases}
\{-1\},&r>0,\\[2mm]
[-1,1],&r=0,\\[2mm]
\{1\},&r<0.
\end{cases}
\]
By \eqref{point}, $-p$ belongs to this concave superdifferential, so
we get, whenever \(|g|>1+\delta\) and \(f\neq0\), $p=\operatorname{sgn}(f)$.
Since \(f_M\) has the same sign as \(f\),
\[
pf_M=|f_M|
\qquad\text{on }\{|g|>1+\delta\}.
\]
Multiplying by \(\chi(g)\ge0\), $p\,f_M\chi(g)=|f_M|\chi(g)$.
Thus
\begin{equation}
\label{pvM}
pv_M=|v_M|.
\end{equation}
On \(\{|g|\le 1+\delta\}\), both sides vanish because \(\chi(g)=0\). If \(f=0\), then \(v_M=0\), so the identity also holds there.

Recall
\[
\eta=\partial_sU_\delta(f,g).
\]
As \(U_\delta(r,s)=a-|r|\) for \(|s|>1+\delta\) we see that $\eta=0$ where $|g|>1+\delta$.
Consequently, everywhere
\[
\eta\chi(g)=0,\qquad\eta\chi'(g)=0.
\]
The latter holds because \(\chi\) is identically zero on \([-(1+\delta),1+\delta]\), so \(\chi'\) can be nonzero only in the exterior region where \(\eta=0\).
Formally, for smooth \(f,g\),
\[
\partial_xv_M
=
\1_{\{|f|<M\}}\,\partial_xf\,\chi(g)
+
f_M\,\chi'(g)\,\partial_xg.
\]
For $W^{1/2,2}$ data, the pairings are between $W^{-1/2,2}$ and
$W^{1/2,2}$. Choose real smooth $f_n,g_n$ converging to $f,g$ in
$W^{1/2,2}$, and set
\[
 \eta_n=\partial_sU_\delta(f_n,g_n),\qquad
 f_{M,n}=\max(-M,\min(f_n,M)),\qquad
 v_{M,n}=f_{M,n}\chi(g_n).
\]
The identities $\eta_n\chi(g_n)=\eta_n\chi'(g_n)=0$ hold pointwise.
The ordinary almost-everywhere chain rule therefore gives
$\eta_n\partial_xv_{M,n}=0$.
By Lemmas~\ref{fp-regularization} and~\ref{fp-composition},
$\eta_n\to\eta$ and $v_{M,n}\to v_M$ in $W^{1/2,2}$.
Also $p=-R(\eta-2\varepsilon g)\in W^{1/2,2}$, since $R$ is
bounded on this space.
Integration by parts and passage to the limit give
\begin{equation}
\label{dxvM}
\langle\partial_x\eta,v_M\rangle
=-\langle\eta,\partial_xv_M\rangle=0.
\end{equation}

The identities \eqref{pvM} and \eqref{dxvM} will be used together with the following chain-rule identity.

For the same approximants, the ordinary chain rule gives
\[
 \langle\partial_xg_n,v_{M,n}\rangle
 =\langle\partial_x\Theta(g_n),f_{M,n}\rangle.
\]
Lipschitz composition gives $\Theta(g_n)\to\Theta(g)$ and
$f_{M,n}\to f_M$ in $W^{1/2,2}$. Passing to the limit in the
$W^{-1/2,2}$--$W^{1/2,2}$ pairings yields
\[
 \langle\partial_xg,v_M\rangle
 =\langle\partial_x\Theta(g),f_M\rangle.
\]
We have proved the three identities:

\begin{align}
\label{vM}
pv_M=|v_M|
&\quad\text{from the cusp }a-|r|,        \notag
\\
\langle\partial_x\eta,v_M\rangle=0
&\quad\text{because }\partial_sU_\delta=0\text{ outside the strip},      \notag
\\
\langle\partial_xg,v_M\rangle
=\langle\partial_x\Theta(g),f_M\rangle
&\quad\text{from }\Theta'=\chi.
\end{align}

\paragraph{The Dirichlet-form estimate.}

Let $\mathsf Q_s=e^{-sA}$. Since \(\mathsf Q_s\) is sub-Markov and \(|p|\le1\),
\[
|\mathsf Q_sp|\le \mathsf Q_s|p|\le \mathsf Q_s1\le1.
\]
Using \(pv_M=|v_M|\) from \eqref{vM}
\[
\begin{aligned}
\langle(I-\mathsf Q_s)p,v_M\rangle
&=\int\bigl(|v_M|-(\mathsf Q_sp)v_M\bigr)\,d\mu \ge 0.
\end{aligned}
\]
Indeed, $(\mathsf Q_sp)v_M\le |\mathsf Q_sp|\,|v_M|\le |v_M|$.
Divide by \(s\) and let \(s\downarrow0\). For the Dirichlet form
\[
\mathcal E(u,v)=\langle A^{1/2}u,A^{1/2}v\rangle,
\]
this gives
\[
0\le\mathcal E(p,v_M).
\]
Using now
\[
Ap=\partial_x\eta-2\varepsilon\partial_xg
\]
and the identities in \eqref{vM},
\[
\begin{aligned}
\mathcal E(p,v_M)
=\langle Ap,v_M\rangle
=\langle\partial_x\eta,v_M\rangle
-2\varepsilon\langle\partial_xg,v_M\rangle
=-2\varepsilon
\langle\partial_x\Theta(g),f_M\rangle.
\end{aligned}
\]

Hence
\[
\langle\partial_x\Theta(g),f_M\rangle\le0.
\]
Letting \(M\to\infty\), using $f_M\to f$ in $W^{1/2,2}$ by Lemma~\ref{fp-composition},
\begin{equation*}
\langle\partial_x\Theta(g),f\rangle\le0.
\end{equation*}
On the other hand, \(Ag=-\partial_xf\), so
\[
\begin{aligned}
\mathcal E(g,\Theta(g))
=\langle Ag,\Theta(g)\rangle
=-\langle\partial_xf,\Theta(g)\rangle
=\langle f,\partial_x\Theta(g)\rangle
\le 0.
\end{aligned}
\]
But \(\Theta\) is nondecreasing and one-Lipschitz. Therefore the Markov contraction property gives
\begin{equation}
\label{Theta}
0\le
\mathcal E(\Theta(g),\Theta(g))
\le
\mathcal E(g,\Theta(g)).
\end{equation}
To verify \eqref{Theta}, apply \eqref{theta} with $s=g(x,a), r=g(y,b)$ and integrate in all variables against the symmetric kernel of $e^{-tA}$, with $t>0$. For a sub-Markov kernel, the missing-mass term is
\[
\int(1-e^{-tA}1)g\Theta(g)\,d\mu
\ge\int(1-e^{-tA}1)\Theta(g)^2\,d\mu,
\]
by the second inequality in \eqref{theta}. Considering $I-e^{-tA}$, dividing by $t$, and letting $t\downarrow0$, we obtain \eqref{Theta}.
Combining this with the preceding inequality,
\[
0\le
\mathcal E(\Theta(g),\Theta(g))
\le
\mathcal E(g,\Theta(g))
\le0.
\]
Thus
\begin{equation*}
\mathcal E(\Theta(g),\Theta(g))=0.
\end{equation*}
Since \(A\ge D=|D_x|\) in the quadratic-form sense, this implies that \(\Theta(g)\) is horizontally constant in each state \(a\in E\):
\[
\Theta(g(x,a))\equiv c_a.
\]
If \(c_a>0\), then \(g(x,a)>1+\delta\) for almost every \(x\), contradicting the fact that \(g=Rf\) has zero horizontal mean. Similarly, \(c_a<0\) would imply \(g(x,a)<-(1+\delta)\) almost everywhere. Therefore \(c_a=0\), and hence
\[
\Theta(g)=0.
\]
By the definition of \(\Theta\),

\[
|g|=|Rf|\le 1+\delta\quad\text{almost everywhere}.
\]
This proves \eqref{Rfde}.

\begin{proof}[Proof of Theorem~\ref{fp-finite}]
The confinement estimate \eqref{Rfde} completes the argument for
\eqref{UR} in Subsection~\ref{fp-remove-section}. This gives
\eqref{fp-finite-U}, including the sub-Markov case by the additional
nonnegative terms in the reversible pairing above.
Since $\Bell(r,s)=U(r,s)+|r|$, apply the inequality to $f/\lambda$
and multiply by $\lambda$ to obtain \eqref{fp-finite-Bellman}.
\end{proof}

\subsection{Passage to general transverse spaces and sharpness}
\label{tr-section}

The constant in Theorem~\ref{fp-finite} is independent of the finite
space, its weights and its generator. We now pass to the general
setting by approximation. Recall that $\Bell$ is continuous and
nonnegative, $\Bell(0,0)=0$, and
\[
 \Bell(r,s)=\cD\quad\text{when }|s|\geq1.
\]
Thus almost-everywhere convergence of the input and its transform
allows passage to the limit by Fatou's lemma. We first pass from the
circle to the line, then approximate the transverse space and extend
to $L^1$. The circle has normalized Haar measure and the line has
Lebesgue measure.

\begin{lemma}\label{tr-line}
Let $E$ be a finite transverse space with positive weights, and let
$-L$ generate a symmetric sub-Markov semigroup on $E$.  Then, on either
$X=\T\times E$ or $X=\R\times E$,
\begin{equation}\label{tr-finite-integral}
 \lambda\int_X
 \Bell\!\left(\frac{f}{\lambda},\frac{R_Lf}{\lambda}\right)d\mu
 \leq\|f\|_1,
 \qquad f\in L^1(X;\R)\cap L^2(X;\R),\quad\lambda>0.
\end{equation}
No normalization of the total transverse weight is required.
\end{lemma}

\begin{proof}
The circle assertion follows from Theorem~\ref{fp-finite} by scaling
$f$. Normalize the transverse measure and multiply the inequality
by its total mass.

For the line, first take $f\in C_c^\infty(\R;\R^E)$ and
$\lambda=1$. Periodize $f$ on a circle of length $\ell$, and denote
the resulting function and transform by $f_\ell$ and $R_{L,\ell}$.
The substitution $\theta=2\pi x/\ell$ replaces $L$ in the transform
by $(\ell/(2\pi))^2L$.  Apply the circle result to this generator and
multiply its normalized integral by $\ell$.  For all sufficiently
large $\ell$ this gives
\begin{equation}\label{tr-periodized}
 \int_{-\ell/2}^{\ell/2}\sum_{j\in E}\nu_j
       \Bell(f_{\ell,j},(R_{L,\ell}f_\ell)_j)\,dx
 \leq\|f\|_1.
\end{equation}

Diagonalize the finite-dimensional self-adjoint operator $L$.
On an eigenspace with eigenvalue $s\geq0$, the transform has multiplier
\[
 m_s(\xi)=-i\xi(\xi^2+s)^{-1/2},\qquad m_s(0)=0.
\]
With the convention $\widehat h(\xi)=\int_\R e^{-ix\xi}h(x)\,dx$,
the corresponding Fourier series is
\[
 \frac1\ell\sum_{k\in\mathbb Z}
 m_s(2\pi k/\ell)\widehat h(2\pi k/\ell)e^{2\pi ikx/\ell}.
\]
This converges locally uniformly to the Fourier integral for
$m_s(D)h$. Away from zero, use Riemann sums. When $s=0$, the integral
and sum over $|\xi|<\eta$ are bounded by $C(\eta+\ell^{-1})$,
uniformly on compact $x$-intervals. The rapid decay of $\widehat h$
controls the tails uniformly in $\ell$. There are finitely many
eigenspaces, so $R_{L,\ell}f_\ell\to R_Lf$ locally uniformly.
Fatou's lemma applied to the nonnegative integrands in
\eqref{tr-periodized}, extended by zero outside their intervals,
proves \eqref{tr-finite-integral} for smooth compactly supported $f$.

For general $f\in L^1\cap L^2$, choose smooth compactly supported
approximants converging in both norms.  The $L^2$ contraction of $R_L$
gives convergence of their transforms in $L^2$.  A subsequence
converges almost everywhere simultaneously with its transforms, and
another application of Fatou proves the result.  Scaling gives every
$\lambda>0$.
\end{proof}

For the transverse approximation, we compress $e^{-rL}$, not $L$:
partition indicators need not lie in the domain of $L$.
The compression retains positivity and the sub-Markov
property, and we then let $r\downarrow0$.

\begin{proposition}\label{tr-compression}
Let $(Y,\nu)$ and $L$ satisfy the assumptions of
Theorem~\ref{thm:horizontal}.  On $X=\T\times Y$ or $\R\times Y$,
\begin{equation}\label{tr-general-integral}
 \lambda\int_X
 \Bell\!\left(\frac{f}{\lambda},\frac{R_Lf}{\lambda}\right)d\mu
 \leq\|f\|_1
 \qquad(f\in L^1(X;\R)\cap L^2(X;\R),\ \lambda>0).
\end{equation}
\end{proposition}

\begin{proof}
Choose increasing finite-measure sets $Y_n$ exhausting $Y$ and finite
partitions $\mathcal P_n$ of $Y_n$ into sets of positive measure.
They can be chosen so that the averaging projections
\[
 E_nh=\sum_{F\in\mathcal P_n}
      \left(\frac1{\nu(F)}\int_F h\,d\nu\right)\mathbf1_F
\]
converge strongly to the identity on both $L^1(Y)$ and $L^2(Y)$.
Indeed, standardness and $\sigma$-finiteness give a countable family
of finite-measure sets whose simple functions are dense in both
spaces. Refine the partitions to include successively these sets.
The projections are contractions and eventually fix every simple
function in the resulting dense union. Write $E_n$ also for the
fiberwise projection on $X$; density of finite sums of product
functions gives the same convergence there.

Fix $r>0$, write $V_n=\operatorname{ran}E_n$, and set
\begin{equation*}
 Q_{n,r}=E_ne^{-rL}E_n\big|_{V_n},\qquad
 L_{n,r}=r^{-1}(I_{V_n}-Q_{n,r}).
\end{equation*}
The space $V_n$ is a finite weighted space, with weights $\nu(F)$.
Its operator $Q_{n,r}$ is positive, self-adjoint and sub-Markov.
It is also nonnegative as a Hilbert-space operator, because
$\langle Q_{n,r}h,h\rangle=\langle e^{-rL}h,h\rangle\geq0$.
Consequently $0\leq L_{n,r}\leq r^{-1}I_{V_n}$, and
\begin{equation*}
 e^{-tL_{n,r}}=e^{-t/r}\sum_{m=0}^\infty
              \frac{(t/r)^m}{m!}Q_{n,r}^{m},\qquad
 e^{-tL_{n,r}}\mathbf1_{Y_n}\leq\mathbf1_{Y_n}.
\end{equation*}
Lemma~\ref{tr-line} therefore applies with the same constant for
every $n,r$, including any loss of mass in the compression.

For spectral convergence only, extend $L_{n,r}$ by zero on $V_n^\perp$.
The extended operator is
$r^{-1}(E_n-E_ne^{-rL}E_n)$, and hence converges strongly, as $n\to\infty$,
to the bounded self-adjoint operator
\[
 L_r=r^{-1}(I-e^{-rL}).
\]
All these operators have spectrum in $[0,r^{-1}]$.
For each nonzero horizontal frequency $\xi$, bounded continuous
functional calculus gives strong convergence of
$-i\xi(\xi^2+L_{n,r})^{-1/2}$ to
$-i\xi(\xi^2+L_r)^{-1/2}$.  At frequency zero all multipliers are
defined to be zero.  Their norms are at most one.  Plancherel's theorem
and dominated convergence over the horizontal frequencies therefore
give
\begin{equation}\label{tr-first-spectral-limit}
 R_{L_{n,r}}E_nf\longrightarrow R_{L_r}f
       \quad\text{in }L^2(X)\qquad(n\to\infty).
\end{equation}
On the left, the finite-space transform is extended by zero outside
$Y_n$; this agrees with applying the spectrally extended operator to
$E_nf$.

Next, $r^{-1}(1-e^{-rs})\to s$ for every $s\geq0$.  The joint spectral
calculus for the horizontal derivative and $L$, again with multiplier
bound one, gives
\begin{equation}\label{tr-second-spectral-limit}
 R_{L_r}f\longrightarrow R_Lf
       \quad\text{in }L^2(X)\qquad(r\downarrow0).
\end{equation}
Apply the finite-space inequality to $E_nf$ at fixed $r$.
Because $E_nf\to f$ in $L^1\cap L^2$, a subsequence in
\eqref{tr-first-spectral-limit} converges almost everywhere together
with $E_nf$.  Fatou's lemma gives the inequality with $R_{L_r}f$.
The integrand outside $Y_n$ contributes zero, since $\Bell(0,0)=0$.
Finally use \eqref{tr-second-spectral-limit} and Fatou once more.
Since the projections are $L^1$ contractions, the right side remains
$\|f\|_1$. This proves \eqref{tr-general-integral}.
\end{proof}

To remove the $L^2$ assumption, apply the weak estimate to
differences. The following lemma gives the extension argument.

\begin{lemma}\label{tr-extension}
Let $(X,\mu)$ and $(Z,\omega)$ be $\sigma$-finite measure spaces, and
let $V,W$ be finite-dimensional normed spaces over $\R$ or $\C$.
Let $\mathcal D\subset L^1(X;V)\cap L^2(X;V)$ be a linear subspace
dense in $L^1(X;V)$.  Suppose that the linear operator
$T:\mathcal D\to L^0(Z;W)$ satisfies
\begin{equation}\label{tr-abstract-weak}
 \omega\{\|Tf\|_W>\lambda\}
      \leq \frac C\lambda\|f\|_{L^1(X;V)}
 \qquad(f\in\mathcal D,\ \lambda>0).
\end{equation}
Then $T$ has a unique linear extension to $L^1(X;V)$ that is continuous
from $L^1$ convergence to convergence in measure, and the extension
satisfies the same estimate.  Here convergence in measure means
$\omega\{\|u_j-u\|_W>\varepsilon\}\to0$ for every $\varepsilon>0$,
also when $\omega(Z)=\infty$.

If $T$ also has an $L^2$-continuous realization and $\mathcal D$ is
dense in $L^1\cap L^2$ for the sum of the two norms, the extension
agrees with that realization on $L^1\cap L^2$.

Suppose in addition that $X=Z$, $\mu=\omega$, and that
$\Gamma:V\times W\to[0,\infty)$ is continuous, with
$\Gamma(0,0)=0$.  If
\begin{equation}\label{tr-cost-extension}
 \int_X\Gamma(f,Tf)\,d\mu\leq C_0\|f\|_1
 \qquad(f\in\mathcal D),
\end{equation}
then this inequality also holds for every $f\in L^1(X;V)$.
\end{lemma}

\begin{proof}
Let $f_j\in\mathcal D$ converge to $f$ in $L^1$.
Applying \eqref{tr-abstract-weak} to differences shows that $Tf_j$ is
Cauchy in measure.  Choose a subsequence, still denoted by $f_j$, such
that
\[
 \sum_j2^j\|f_{j+1}-f_j\|_1<\infty.
\]
The measures of the sets
$\{\|Tf_{j+1}-Tf_j\|_W>2^{-j}\}$ have finite sum.
Outside their null limsup, the sequence $Tf_j$ is pointwise Cauchy.
Since $W$ is complete, it converges almost everywhere to a finite
measurable function, which we denote by $Tf$.

This limit is the limit in measure of the original sequence.
Indeed, for fixed $j$ and $\varepsilon>0$, pointwise convergence along
the chosen subsequence and \eqref{tr-abstract-weak} give
\[
 \omega\{\|Tf_j-Tf\|_W>\varepsilon\}
 \leq\liminf_{k\to\infty}
      \omega\{\|Tf_j-Tf_k\|_W>\varepsilon\}
 \leq\frac C\varepsilon\|f_j-f\|_1.
\]
Using two approximating sequences shows that the limit is independent
of the approximation. The same estimate proves continuity; passage
to limits in the identities on $\mathcal D$ gives linearity.
Density then gives uniqueness.

An almost-everywhere convergent subsequence and Fatou's lemma for
the indicators of the open sets $\{\|w\|_W>\lambda\}$ give
\[
 \omega\{\|Tf\|_W>\lambda\}
 \leq\liminf_j\omega\{\|Tf_j\|_W>\lambda\}
 \leq\frac C\lambda\|f\|_1.
\]
If an $L^2$ realization is available, simultaneous $L^1,L^2$
approximation and $L^2$ continuity give the same limit in measure.
For \eqref{tr-cost-extension}, choose a subsequence on which both
$f_j\to f$ and $Tf_j\to Tf$ almost everywhere. Continuity and
nonnegativity of $\Gamma$ allow another application of Fatou's lemma.
\end{proof}

\begin{proposition}\label{tr-upper}
The integral and upper weak-type estimates in
Theorems~\ref{thm:horizontal} and~\ref{thm:horizontal-integral}
hold for every real $f\in L^1(X)$.
\end{proposition}

\begin{proof}
By \eqref{tr-general-integral} and the value of $\Bell$ outside the
strip,
\[
 \lambda\mu\{|R_Lf|>\lambda\}\leq\KD\|f\|_1
       \qquad(f\in L^1\cap L^2).
\]
Apply Lemma~\ref{tr-extension} with $\mathcal D=L^1\cap L^2$.
Its last assertion, with
$\Gamma(r,s)=\lambda\Bell(r/\lambda,s/\lambda)$, preserves the
integral inequality for each $\lambda>0$.
The pointwise bounds
$\Phi(|s|)\leq\Bell(r,s)$ and
\[
 \Phi(t)\geq t^2/2\quad(0\leq t<1),\qquad
 \Phi(t)=\cD\quad(t\geq1),
\]
then give \eqref{eq:main-distribution} and
\eqref{eq:main-truncated}.
\end{proof}

For sharpness, use Davis's example at high horizontal frequencies.
At each fixed transverse spectral value the multiplier then tends to
that of the Hilbert transform. No nonzero vector in $\ker L$ is needed.

\begin{proposition}\label{tr-sharpness}
For every nonzero transverse space in Theorem~\ref{thm:horizontal},
on either the circle or the line product,
\[
 \|R_L\|_{L^1(X;\R)\to L^{1,\infty}(X)}\geq\KD.
\]
Consequently the upper weak-type constant in that theorem is exact,
even for semigroups with killing.
\end{proposition}

\begin{proof}
Use Davis's circle pair \cite{Davis1974}, with the horizontal Hilbert
transform convention fixed above:
\begin{equation*}
 g_*(x)=\frac2\pi\log\left|\cot\frac x2\right|,\qquad
 \Hilb g_*=h_*:=\operatorname{sgn}(\sin x),\qquad
 \|g_*\|_{L^1(\T)}=\cD,\qquad g_*\in L^2(\T).
\end{equation*}
The identities also follow from the Fourier expansion
$g_*=(4/\pi)\sum_{k\geq1,\ k\text{ odd}}\cos(kx)/k$ in $L^2$.
Choose $F\subset Y$ with $0<\nu(F)<\infty$ and put $w=\mathbf1_F$.
For each nonzero integer $k$, the spectral theorem gives
\begin{equation}\label{tr-circle-frequency}
 -ikN\bigl((kN)^2+L\bigr)^{-1/2}w
       \longrightarrow-i\operatorname{sgn}(k)w
       \quad\text{in }L^2(Y)\qquad(N\to\infty).
\end{equation}
The scalar multipliers converge pointwise on $[0,\infty)$ and are
bounded by one; no assumption on the semigroup at large times is
needed. Approximate $g_*$ in $L^2(\T)$ by mean-zero trigonometric
polynomials and apply \eqref{tr-circle-frequency} term by term.
The $L^2$ contractions of $R_L$ and $\Hilb$ control the tails uniformly
in integer $N$. Thus
\[
 R_L\bigl(g_*(Nx)w\bigr)-h_*(Nx)w\longrightarrow0
       \quad\text{in }L^2(\T\times Y).
\]
The input norm is $\cD\nu(F)$, whereas $|h_*(Nx)w|=1$ almost
everywhere on $\T\times F$.  For $0<\lambda<1$, Chebyshev's
inequality for the error therefore gives
\[
 \liminf_{N\to\infty}
 \mu\{|R_L(g_*(Nx)w)|>\lambda\}\geq\nu(F).
\]
First let $N\to\infty$ and then $\lambda\uparrow1$ to obtain the
circle lower bound $1/\cD=\KD$.

On the line, cut off the periodic functions by
$\varphi=\mathbf1_{[0,2\pi]}$.
For every nonzero integer $k$,
\begin{equation}\label{tr-line-frequency}
 \bigl\|R_L(e^{ikNx}\varphi w)
       +i\operatorname{sgn}(k)e^{ikNx}\varphi w\bigr\|_2
       \longrightarrow0.
\end{equation}
After translating the horizontal Fourier variable, the relevant
multiplier is
\[
 -i(\xi+kN)\bigl((\xi+kN)^2+s\bigr)^{-1/2},
\]
which converges to $-i\operatorname{sgn}(k)$ for every $\xi\in\R$
and $s\geq0$, and is bounded in modulus by one.
Dominated convergence against the joint spectral measure of
$\varphi w$ proves \eqref{tr-line-frequency}.
Fourier approximation of $g_*$ now gives
\[
 R_L\bigl(\varphi g_*(Nx)w\bigr)
       -\varphi h_*(Nx)w\longrightarrow0\quad\text{in }L^2(\R\times Y).
\]
The approximation tails are uniform in integer $N$, since
$[0,2\pi]$ contains exactly $N$ periods of $g_*(Nx)$.
The input norm is $2\pi\cD\nu(F)$ and the target has modulus one
on a set of measure $2\pi\nu(F)$.  The preceding level-set argument
again gives the lower bound $\KD$.
\end{proof}

\section{Applications}
\label{ex-section}

Theorem~\ref{thm:horizontal} applies whenever the transverse
semigroup is symmetric and sub-Markov; no heat-kernel bound,
compactness or conservation of mass is needed. Here are several
examples, some of which also occur in
\cite[Corollary 4.1]{CicconeEtAl2026}. In all cases the numerator
$-\partial_x$ acts on an additional line or circle.

\begin{corollary}\label{ex-generators}
Let $d\geq1$.  All conclusions of Theorems~\ref{thm:horizontal}
and~\ref{thm:horizontal-integral} hold
on $X=\R\times Y$ and $X=\T\times Y$ for
\[
 R_L=-\partial_x(-\partial_x^2+L)^{-1/2}
\]
with any of the following transverse generators and measures.

\begin{enumerate}[label=(\alph*)]
\item The Euclidean Laplacian $L=-\Delta_y$ on
      $L^2(\R^d,dy)$.
\item The periodic Laplacian $L=-\Delta_{\T^d}$ on
      $L^2(\T^d)$, with normalized Haar measure.
\item The Dunkl Laplacian $L=-\Delta_\kappa$ on
      $L^2(\R^d,w_\kappa(y)\,dy)$, where
      \[
      w_\kappa(y)=\prod_{\alpha\in\mathcal R_+}
                    |\langle\alpha,y\rangle|^{2\kappa(\alpha)}.
      \]
      Here $\mathcal R$ is a root system, $\mathcal R_+$ a positive
      subsystem, and $\kappa$ a nonnegative multiplicity invariant
      under the associated reflection group.
\item The Ornstein--Uhlenbeck generator
      \[
      L=-\Delta_y+y\cdot\nabla_y
      \quad\text{on }L^2(\R^d,\gamma_d),\qquad
      d\gamma_d(y)=(2\pi)^{-d/2}e^{-|y|^2/2}\,dy.
      \]
\item The Laguerre polynomial generator
      \[
      L=-\sum_{j=1}^d
         \bigl[y_j\partial_{y_j}^2+
                  (\alpha_j+1-y_j)\partial_{y_j}\bigr]
      \]
      on $L^2((0,\infty)^d,\mu_\alpha)$, where
      \[
      d\mu_\alpha(y)=\prod_{j=1}^d
             \frac{y_j^{\alpha_j}e^{-y_j}}{\Gamma(\alpha_j+1)}\,dy_j,
      \qquad\alpha_j>-1.
      \]
\item The Jacobi polynomial generator
      \[
      L=-\sum_{j=1}^d\bigl[(1-y_j^2)\partial_{y_j}^2+
      \{\beta_j-\alpha_j-(\alpha_j+\beta_j+2)y_j\}
                                      \partial_{y_j}\bigr]
      \]
      on $L^2((-1,1)^d,\mu_{\alpha,\beta})$, where
      \[
      d\mu_{\alpha,\beta}(y)=\prod_{j=1}^d
      \frac{(1-y_j)^{\alpha_j}(1+y_j)^{\beta_j}}{Z_j}\,dy_j,
      \quad Z_j=\int_{-1}^1(1-t)^{\alpha_j}(1+t)^{\beta_j}\,dt,
      \]
      and $\alpha_j,\beta_j>-1$.
\item The Schr\"odinger operator $L=-\Delta_y+V(y)$ on
      $L^2(\R^d,dy)$, where $V\geq0$ and $V\in L^1_{\rm loc}(\R^d)$,
      defined by the closed form
      \[
      \mathfrak q_V[v]=\int_{\R^d}(|\nabla v|^2+V|v|^2)\,dy,
      \qquad
      \operatorname{Dom}(\mathfrak q_V)
             =H^1(\R^d)\cap L^2(\R^d,V\,dy).
      \]
      This includes the unshifted harmonic oscillator
      $-\Delta_y+|y|^2$ and the massive Laplacian
      $-\Delta_y+m^2$, $m\geq0$.
\item The nonnegative Laplace--Beltrami operator $L=-\Delta_M$ on
      $L^2(M,d\operatorname{vol}_M)$, where $M$ is a complete connected
      Riemannian manifold of positive dimension without boundary.
      No curvature bound or stochastic-completeness assumption is
      required.
\end{enumerate}
For the two polynomial examples, use the self-adjoint realizations
diagonalized by the orthogonal polynomial bases.  Their
one-dimensional eigenvalues are $k$ for Laguerre and
$k(k+\alpha+\beta+1)$ for Jacobi, $k\geq0$; the product generators
are sums of the one-dimensional realizations.

In particular, in every listed case,
\begin{equation*}
 \|R_L\|_{L^1(X;\R)\to L^{1,\infty}(X)}
       =\KD=\frac{\pi^2}{8\Gcat}.
\end{equation*}
The same conclusions hold with $L$ replaced by $\psi(L)$ for any
nonnegative Bernstein function $\psi$, including bounded functions
and functions with $\psi(0)>0$.
\end{corollary}

\begin{proof}
The Euclidean and periodic heat semigroups are symmetric and
conservative; see \cite{Stein1970}.
For nonnegative multiplicities, the Dunkl heat semigroup is positive,
symmetric and conservative
\cite[Section 4.2]{Rosler2003}; see also \cite{RoslerVoit1998}.
Mehler's formula gives the same properties for the
Ornstein--Uhlenbeck semigroup
\cite[Chapter 2]{UrbinaRomero2019}.

The Laguerre and Jacobi polynomial semigroups are positive symmetric
contractions \cite[Section 2]{NowakStempak2012},
\cite[Section 2]{NowakSjogren2013}.  The constant polynomial is an
eigenfunction of eigenvalue zero, so both preserve constants.
The full parameter ranges stated above concern these polynomial
Markov realizations, rather than other conjugated or shifted
Laguerre and Jacobi operators.

The Schr\"odinger form is densely defined, closed and Markovian.
Its associated semigroup is therefore symmetric and sub-Markov
\cite[Theorem 1.4.1]{FukushimaOshimaTakeda2011}; see also
\cite[Section 2]{KucharskiWrobel2026} for positivity and domination
by the heat semigroup.  The heat semigroup on a complete manifold is
symmetric and sub-Markov even when it loses mass at infinity
\cite[Theorem 3.3 and equation (3.11)]{Grigoryan2006}.
Thus Theorem~\ref{thm:horizontal} applies in every case.

For the subordination assertion, use Bochner's formula
\begin{equation*}
 e^{-t\psi(L)}=\int_{[0,\infty)}e^{-sL}\,d\eta_t(s),
 \qquad \eta_t([0,\infty))=e^{-t\psi(0)}\leq1,
\end{equation*}
where $\eta_t$ is a positive subordination measure
\cite[Theorem 5.2 and Chapter 13]{SchillingSongVondracek2012}.
Positivity, symmetry and sub-Markov contractivity are preserved.
Strong continuity follows from the spectral theorem.
Thus the theorem applies to $\psi(L)$. When $\psi(0)>0$, the
subordination measure has mass less than one, corresponding to
additional killing.
\end{proof}

The subordination clause includes $L^\alpha$ for $0<\alpha\leq1$,
$m^2+L$ for $m\geq0$, and the bounded generator $I-e^{-L}$.

In the Schr\"odinger example the conclusion concerns
\[
 -\partial_x\bigl(-\partial_x^2-\Delta_y+V(y)\bigr)^{-1/2},
\]
where the potential is independent of $x$.  On $\R^{1+d}$ this is,
up to sign, the horizontal component of the Schr\"odinger gradient
Riesz transform; compare its $L^p$ theory in
\cite[Section 1]{AuscherBenAli2007}.
For the manifold example, the product metric on
$N=\R\times M$ or $\T\times M$ gives
\[
 -\Delta_N=-\partial_x^2-\Delta_M,
 \qquad R_L=-\partial_x(-\Delta_N)^{-1/2}.
\]
Thus this is the horizontal component of the Riemannian Riesz
transform on $N$, with our sign convention; for the manifold $L^p$
theory under geometric and heat-kernel hypotheses, see
\cite[Sections 1.1--1.2]{AuscherCoulhonDuongHofmann2004}.

In particular, these examples do not assert sharp norms for
transforms whose numerators differentiate in the transverse
variables, such as intrinsic Gaussian, Schr\"odinger or manifold
Riesz transforms.


\begin{thebibliography}{99}
\small

\bibitem{AuscherBenAli2007}
P.~Auscher and B.~Ben Ali,
Maximal inequalities and Riesz transform estimates on \(L^p\) spaces
for Schr\"odinger operators with nonnegative potentials,
\emph{Ann. Inst. Fourier (Grenoble)} \textbf{57} (2007), no.~6, 1975--2013.
doi:\,\href{https://doi.org/10.5802/aif.2320}{10.5802/aif.2320}.

\bibitem{AuscherCoulhonDuongHofmann2004}
P.~Auscher, T.~Coulhon, X.~T.~Duong, and S.~Hofmann,
Riesz transform on manifolds and heat kernel regularity,
\emph{Ann. Sci. \'Ecole Norm. Sup.} (4) \textbf{37} (2004), no.~6, 911--957.
doi:\,\href{https://doi.org/10.1016/j.ansens.2004.10.003}{\nolinkurl{10.1016/j.ansens.2004.10.003}}.

\bibitem{Baernstein1978}
A.~Baernstein II,
Some sharp inequalities for conjugate functions,
\emph{Indiana Univ. Math. J.} \textbf{27} (1978), no.~5, 833--852.
doi:\,\href{https://doi.org/10.1512/iumj.1978.27.27055}{10.1512/iumj.1978.27.27055}.

\bibitem{Banuelos2010}
R.~Ba\~nuelos,
The foundational inequalities of D.~L.~Burkholder and some of their ramifications,
\emph{Illinois J. Math.} \textbf{54} (2010), no.~3, 789--868 (2012).
doi:\,\href{https://doi.org/10.1215/ijm/1336049979}{10.1215/ijm/1336049979}.

\bibitem{BanuelosBaudoinChenSire2021}
R.~Ba\~nuelos, F.~Baudoin, L.~Chen, and Y.~Sire,
Multiplier theorems via martingale transforms,
\emph{J. Funct. Anal.} \textbf{281} (2021), no.~9, Paper No.~109188, 37 pp.
doi:\,\href{https://doi.org/10.1016/j.jfa.2021.109188}{10.1016/j.jfa.2021.109188}.

\bibitem{BanuelosOsekowski2015}
R.~Ba\~nuelos and A.~Os\k{e}kowski,
Sharp martingale inequalities and applications to Riesz transforms
on manifolds, Lie groups and Gauss space,
\emph{J. Funct. Anal.} \textbf{269} (2015), no.~6, 1652--1713.
doi:\,\href{https://doi.org/10.1016/j.jfa.2015.06.015}{10.1016/j.jfa.2015.06.015}.

\bibitem{BanuelosWang1995}
R.~Ba\~nuelos and G.~Wang,
Sharp inequalities for martingales with applications to the Beurling--Ahlfors
and Riesz transforms,
\emph{Duke Math. J.} \textbf{80} (1995), no.~3, 575--600.
doi:\,\href{https://doi.org/10.1215/S0012-7094-95-08020-X}{10.1215/S0012-7094-95-08020-X}.

\bibitem{Burkholder1984}
D.~L.~Burkholder,
Boundary value problems and sharp inequalities for martingale transforms,
\emph{Ann. Probab.} \textbf{12} (1984), no.~3, 647--702.
doi:\,\href{https://doi.org/10.1214/aop/1176993220}{10.1214/aop/1176993220}.

\bibitem{CZ}
A.~P.~Calder\'on and A.~Zygmund,
On the existence of certain singular integrals,
\emph{Acta Math.} \textbf{88} (1952), 85--139.
doi:\,\href{https://doi.org/10.1007/BF02392130}{10.1007/BF02392130}.

\bibitem{CarbonaroDragicevic2013}
A.~Carbonaro and O.~Dragi\v{c}evi\'c,
Bellman function and linear dimension-free estimates in a theorem of Bakry,
\emph{J. Funct. Anal.} \textbf{265} (2013), no.~7, 1085--1104.
doi:\,\href{https://doi.org/10.1016/j.jfa.2013.05.031}{10.1016/j.jfa.2013.05.031}.

\bibitem{CicconeEtAl2026}
V.~Ciccone, M.~Hormozi, J.~Madrid, J.~Niksi\'nski, and B.~Wr\'obel,
On \(L^p\) norms of maximal operators,
preprint, 2026, \href{https://arxiv.org/abs/2610.10247}{arXiv:2610.10247}.

\bibitem{Davis1974}
B.~Davis,
On the weak type \((1,1)\) inequality for conjugate functions,
\emph{Proc. Amer. Math. Soc.} \textbf{44} (1974), no.~2, 307--311.
doi:\,\href{https://doi.org/10.1090/S0002-9939-1974-0348381-6}{10.1090/S0002-9939-1974-0348381-6}.

\bibitem{DomelevoPetermichlSkreb2025}
K.~Domelevo, S.~Petermichl, and K.~A.~\v{S}kreb,
Continuous sparse domination and dimensionless weighted estimates
for the Bakry--Riesz vector,
\emph{J. Reine Angew. Math.} \textbf{824} (2025), 137--166.
doi:\,\href{https://doi.org/10.1515/crelle-2025-0024}{10.1515/crelle-2025-0024}.

\bibitem{FukushimaOshimaTakeda2011}
M.~Fukushima, Y.~Oshima, and M.~Takeda,
\emph{Dirichlet Forms and Symmetric Markov Processes}, 2nd ed.,
De Gruyter, Berlin, 2011.
doi:\,\href{https://doi.org/10.1515/9783110218091}{10.1515/9783110218091}.

\bibitem{GohbergKrupnik1968}
I.~Ts.~Gohberg and N.~Ya.~Krupnik,
Norm of the Hilbert transformation in the \(L_p\) space,
\emph{Funct. Anal. Appl.} \textbf{2} (1968), no.~2, 180--181.
doi:\,\href{https://doi.org/10.1007/BF01075955}{10.1007/BF01075955}.

\bibitem{Grigoryan2006}
A.~Grigor'yan,
Heat kernels on weighted manifolds and applications,
in \emph{The Ubiquitous Heat Kernel},
Contemporary Mathematics, vol.~398,
American Mathematical Society, Providence, RI, 2006, pp.~93--191.
doi:\,\href{https://doi.org/10.1090/conm/398/07486}{10.1090/conm/398/07486}.

\bibitem{GundyVaropoulos1979}
R.~F.~Gundy and N.~Th.~Varopoulos,
Les transformations de Riesz et les int\'egrales stochastiques,
\emph{C. R. Acad. Sci. Paris S\'er. A-B} \textbf{289} (1979), no.~1, A13--A16.

\bibitem{IwMa}
T.~Iwaniec and G.~Martin,
Riesz transforms and related singular integrals,
\emph{J. Reine Angew. Math.} \textbf{473} (1996), 25--57.
doi:\,\href{https://doi.org/10.1515/crll.1995.473.25}{10.1515/crll.1995.473.25}.

\bibitem{JanakiramanWeak2004}
P.~Janakiraman,
Best weak-type $(p,p)$ constants, $1\le p\le2$, for orthogonal
harmonic functions and martingales,
\emph{Illinois J. Math.} \textbf{48} (2004), no.~3, 909--921.
doi:\,\href{https://doi.org/10.1215/ijm/1258131059}{10.1215/ijm/1258131059}.

\bibitem{Janakiraman2004}
P.~Janakiraman,
Weak-type estimates for singular integrals and the Riesz transform,
\emph{Indiana Univ. Math. J.} \textbf{53} (2004), no.~2, 533--556.
doi:\,\href{https://doi.org/10.1512/IUMJ.2004.53.2372}{10.1512/IUMJ.2004.53.2372}.

\bibitem{Kol}
A.~Kolmogoroff,
Sur les fonctions harmoniques conjugu\'ees et les s\'eries de Fourier,
\emph{Fund. Math.} \textbf{7} (1925), 24--29.
doi:\,\href{https://doi.org/10.4064/fm-7-1-24-29}{10.4064/fm-7-1-24-29}.

\bibitem{KucharskiWrobel2026}
M.~Kucharski and B.~Wr\'obel,
On \(L^p\) estimates for positivity-preserving Riesz transforms related
to Schr\"odinger operators,
\emph{Ann. Inst. Fourier (Grenoble)} \textbf{76} (2026), no.~4, 1593--1634.
doi:\,\href{https://doi.org/10.5802/aif.3744}{10.5802/aif.3744}.

\bibitem{NazarovTreilVolberg1999}
F.~Nazarov, S.~Treil, and A.~Volberg,
The Bellman functions and two-weight inequalities for Haar multipliers,
\emph{J. Amer. Math. Soc.} \textbf{12} (1999), no.~4, 909--928.
doi:\,\href{https://doi.org/10.1090/S0894-0347-99-00310-0}{10.1090/S0894-0347-99-00310-0}.

\bibitem{NTV}
F.~Nazarov, S.~Treil, and A.~Volberg,
Bellman function in stochastic control and harmonic analysis,
in \emph{Systems, Approximation, Singular Integral Operators, and Related Topics},
Operator Theory: Advances and Applications, vol.~129,
Birkh\"auser, Basel, 2001, pp.~393--423.
doi:\,\href{https://doi.org/10.1007/978-3-0348-8362-7_16}{10.1007/978-3-0348-8362-7\_16}.


\bibitem{NowakSjogren2013}
A.~Nowak and P.~Sj\"ogren,
Sharp estimates of the Jacobi heat kernel,
\emph{Studia Math.} \textbf{218} (2013), no.~3, 219--244.
doi:\,\href{https://doi.org/10.4064/sm218-3-2}{10.4064/sm218-3-2}.

\bibitem{NowakStempak2012}
A.~Nowak and K.~Stempak,
On \(L^p\)-contractivity of Laguerre semigroups,
\emph{Illinois J. Math.} \textbf{56} (2012), no.~2, 433--452.
doi:\,\href{https://doi.org/10.1215/ijm/1385129958}{10.1215/ijm/1385129958}.

\bibitem{OsekowskiLog2012}
A.~Os\k{e}kowski,
Sharp logarithmic inequalities for Riesz transforms,
\emph{J. Funct. Anal.} \textbf{263} (2012), no.~1, 89--108.
doi:\,\href{https://doi.org/10.1016/j.jfa.2012.04.007}{10.1016/j.jfa.2012.04.007}.

\bibitem{Osekowski2012}
A.~Os\k{e}kowski,
\emph{Sharp Martingale and Semimartingale Inequalities},
Monografie Matematyczne, vol.~72,
Birkh\"auser/Springer Basel AG, Basel, 2012.
doi:\,\href{https://doi.org/10.1007/978-3-0348-0370-0}{10.1007/978-3-0348-0370-0}.

\bibitem{OsekowskiSurvey2013}
A.~Os\k{e}kowski,
Survey article: Bellman function method and sharp inequalities for martingales,
\emph{Rocky Mountain J. Math.} \textbf{43} (2013), no.~6, 1759--1823.
doi:\,\href{https://doi.org/10.1216/RMJ-2013-43-6-1759}{10.1216/RMJ-2013-43-6-1759}.

\bibitem{OSS2026}
Y.~Ouyang, D.~Spector, and C.~B.~Stockdale,
A dimension-free weak-type \((1,1)\) bound for the vector Riesz transform
on \(\R^n\),
preprint, 2026,
\href{https://arxiv.org/abs/2608.18068}{arXiv:2608.18068}.

\bibitem{Pichorides1972}
S.~K.~Pichorides,
On the best values of the constants in the theorems of M.~Riesz,
Zygmund and Kolmogorov,
\emph{Studia Math.} \textbf{44} (1972), no.~2, 165--179.
doi:\,\href{https://doi.org/10.4064/sm-44-2-165-179}{10.4064/sm-44-2-165-179}.

\bibitem{Riesz}
M.~Riesz,
Sur les fonctions conjugu\'ees,
\emph{Math. Z.} \textbf{27} (1928), no.~1, 218--244.
doi:\,\href{https://doi.org/10.1007/BF01171098}{\nolinkurl{10.1007/BF01171098}}.

\bibitem{Rosler2003}
M.~R\"osler,
Dunkl operators: theory and applications,
in \emph{Orthogonal Polynomials and Special Functions},
Lecture Notes in Mathematics, vol.~1817,
Springer, Berlin, 2003, pp.~93--135.
doi:\,\href{https://doi.org/10.1007/3-540-44945-0_3}{10.1007/3-540-44945-0\_3}.

\bibitem{RoslerVoit1998}
M.~R\"osler and M.~Voit,
Markov processes related with Dunkl operators,
\emph{Adv. Appl. Math.} \textbf{21} (1998), no.~4, 575--643.
doi:\,\href{https://doi.org/10.1006/aama.1998.0609}{10.1006/aama.1998.0609}.

\bibitem{SchillingSongVondracek2012}
R.~L.~Schilling, R.~Song, and Z.~Vondra\v{c}ek,
\emph{Bernstein Functions: Theory and Applications},
2nd revised and extended ed.,
De Gruyter Studies in Mathematics, vol.~37,
De Gruyter, Berlin--Boston, 2012.
doi:\,\href{https://doi.org/10.1515/9783110269338}{10.1515/9783110269338}.

\bibitem{SpectorStockdale2021}
D.~Spector and C.~B.~Stockdale,
On the dimensional weak-type \((1,1)\) bound for Riesz transforms,
\emph{Commun. Contemp. Math.} \textbf{23} (2021), no.~7,
2050072, 19 pp.
doi:\,\href{https://doi.org/10.1142/S0219199720500728}{\nolinkurl{10.1142/S0219199720500728}}.

\bibitem{SteinICM1986}
E.~M.~Stein,
Problems in harmonic analysis related to curvature and oscillatory integrals,
in \emph{Proceedings of the International Congress of Mathematicians,
Berkeley, California, 1986}, vol.~1,
American Mathematical Society, Providence, RI, 1987, pp.~196--221.

\bibitem{Stein1970}
E.~M.~Stein,
\emph{Topics in Harmonic Analysis Related to the Littlewood--Paley Theory},
Annals of Mathematics Studies, vol.~63,
Princeton University Press, Princeton, NJ, 1970.

\bibitem{UrbinaRomero2019}
W.~Urbina-Romero,
\emph{Gaussian Harmonic Analysis},
Springer Monographs in Mathematics, Springer, Cham, 2019.
doi:\,\href{https://doi.org/10.1007/978-3-030-05597-4}{10.1007/978-3-030-05597-4}.

\bibitem{VasyuninVolberg2020}
V.~Vasyunin and A.~Volberg,
\emph{The Bellman Function Technique in Harmonic Analysis},
Cambridge Studies in Advanced Mathematics, vol.~186,
Cambridge University Press, Cambridge, 2020.
doi:\,\href{https://doi.org/10.1017/9781108764469}{10.1017/9781108764469}.
\bibitem{VolbergNazarov2004}
A.~Volberg and F.~Nazarov,
Heating of the Ahlfors--Beurling operator, and estimates of its norm,
\emph{St. Petersburg Math. J.} \textbf{15} (2004), no.~4, 563--573.
doi:\,\href{https://doi.org/10.1090/S1061-0022-04-00822-2}{10.1090/S1061-0022-04-00822-2}.

\end{thebibliography}
\end{document}